\documentclass[a4paper]{amsart}
\usepackage{color}
\usepackage[latin1]{inputenc}
\usepackage[T1]{fontenc}
\usepackage{amsfonts}
\usepackage{amssymb}
\usepackage{amsmath}
\usepackage{amsthm}
\usepackage{stmaryrd}
\usepackage{eufrak}
\usepackage{graphicx,type1cm,eso-pic,color}
\usepackage{pstricks} 
\usepackage{float}
\newtheorem{theorem}{Theorem}[section]
\newtheorem{lemma}[theorem]{Lemma}
\newtheorem{proposition}[theorem]{Proposition}
\newtheorem{corollary}[theorem]{Corollary}
\newtheorem{definition}[theorem]{Definition}
\newtheorem{assumption}[theorem]{Assumption}
\newtheorem{remark}[theorem]{Remark}
\begin{document}
\setlength\arraycolsep{2pt}
\title{Nonparametric Estimation for I.I.D. Paths of Fractional SDE}
\author{Fabienne COMTE*}
\address{*Laboratoire MAP5 UMR CNRS 8145, Universit\'e de Paris, Paris, France}
\email{fabienne.comte@parisdescartes.fr}
\author{Nicolas MARIE$^\dag$}
\address{$^\dag$Laboratoire Modal'X, Universit\'e Paris Nanterre, Nanterre, France}
\email{nmarie@parisnanterre.fr}
\address{$^\dag$ESME Sudria, Paris, France}
\email{nicolas.marie@esme.fr}
\date{}
\maketitle
\noindent
%


%
\begin{abstract}
This paper deals with nonparametric estimators of the drift function $b$ computed from independent continuous observations, on a compact time interval, of the solution of a stochastic differential equation driven by the fractional Brownian motion (fSDE). First, a risk bound is established on a Skorokhod's integral based least squares oracle $\widehat b$ of $b$. Thanks to the relationship between the solution of the fSDE and its derivative with respect to the initial condition, a risk bound is deduced on a calculable approximation of $\widehat b$. Another bound is directly established on an estimator of $b'$ for comparison. The consistency and rates of convergence are established for these estimators in the case of the compactly supported trigonometric basis or the $\mathbb R$-supported Hermite basis.
\end{abstract}


\noindent {\bf MS classification 2020.} 62M09 - 62G08 - 60G22 - 60H07  
%
\tableofcontents
%


%
\section{Introduction}
Consider the stochastic differential equation
\begin{equation}\label{main_equation}
X(t) = x_0 +\int_{0}^{t}b(X(s))ds +\sigma B(t)
\textrm{ $;$ }t\in [0,T],
\end{equation}
where $\sigma,T > 0$, $B$ is a fractional Brownian motion of Hurst index $H\in ]1/2,1[$, $b :\mathbb R\rightarrow\mathbb R$ is a continuous map and $x_0\in\mathbb R^*$.
\\
\\
In this work, we assume that we observe $N$ i.i.d. paths of the solution $X$ of Equation (\ref{main_equation}). For instance, this situation may occur in pharmacokinetics when a group of patients can be monitored: for each patient, a bolus of drug is injected and the "path" of its diffusion in the body can be observed (with a delay) (see details in Subsection \ref{subsubsection_application_context}). Our aim is to propose and study nonparametric estimators of the drift function $b$ based on these observations. This problem is related to functional data analysis, and more specifically, there are various recent contributions about i.i.d. parametric models of (non fractional) stochastic differential equations with mixed effects (see, e.g., Ditlevsen and De Gaetano \cite{DDG05}, Overgaard et al. \cite{OJTM05}, Picchini, De Gaetano and Ditlevsen \cite{PDGD10}, Picchini and Ditlevsen \cite{PD11}, Comte, Genon-Catalot and Samson \cite{CGCS13}, Delattre and Lavielle \cite{DL13}, Delattre, Genon-Catalot and Samson \cite{DGCS13}, Dion and Genon-Catalot \cite{DGC16}, Delattre, Genon-Catalot and Lar\'edo \cite{DGCL18}). Also, i.i.d. samples of stochastic differential equations have been  recently considered in the framework of multiclass classification of diffusions (see Denis, Dion and Martinez \cite{DDM20}). The need of flexibility to deal with the information contained in functional data analysis make it interesting to use a nonparametric approach.
\\
\\
Along the last two decades, many authors studied statistical inference from observations drawn from stochastic differential equations driven by fractional Brownian motion, considering the observation of one path, either in continuous time, or in discrete time with fixed or small  step size.
\\
Most references on the estimation of the trend component in Equation (\ref{main_equation}) deal with parametric estimators. Let us start by papers considering continuous time observations. In Kleptsyna and Le Breton \cite{KB01} and Hu and Nualart \cite{HN10}, strategies to estimate the trend component in Langevin's equation are studied. Kleptsyna and Le Breton \cite{KB01} provide a maximum likelihood estimator, where the stochastic integral with respect to the solution of Equation (\ref{main_equation}) returns to an It\^o integral. In \cite{TV07}, Tudor and Viens extend this estimator to equations with a drift function depending linearly on the unknown parameter. Hu and Nualart \cite{HN10} provide a least squares oracle, not an estimator, because the stochastic integral with respect to the solution of Equation (\ref{main_equation}) is taken in the sense of Skorokhod and is not computable. In \cite{HNZ18}, Hu, Nualart and Zhou extend this oracle to equations with a drift function depending linearly on the unknown parameter. Finally, in \cite{MRF20}, Marie and Raynaud de Fitte extend this oracle to non-homogeneous semi-linear equations with almost periodic coefficients.
\\
Now, considering discrete time observations, still in the parametric context, Tindel and Neuenkirch \cite{NT14} study a least squares-type estimator defined by an objective function, tailor-maid with respect to the main result of Tudor and Viens \cite{TV09} on the rate of convergence of the quadratic variation of the fractional Brownian motion. In \cite{PTV20}, Panloup, Tindel and Varvenne extend the results of \cite{NT14} under much more flexible conditions. In \cite{CT13}, Chronopoulou and Tindel provide a likelihood based numerical procedure to estimate a parameter involved in both the drift and the volatility functions in a stochastic differential equation with multiplicative fractional noise.
\\
About  nonparametric methods for the estimation of the function $b$ in Equation (\ref{main_equation}), there are only few references. Saussereau \cite{SAUSSEREAU14} and Comte and Marie \cite{CM19} study the consistency of Nadaraya-Watson type estimators of the drift function $b$ in Equation (\ref{main_equation}). In \cite{MP11}, Mishra and Prakasa Rao established the consistency and a rate of convergence of a nonparametric estimator of the whole trend of the solution to Equation (\ref{main_equation}) extending that of Kutoyants \cite{KUTOYANTS94}. Marie \cite{MARIE19} deals with the same estimator but for reflected fractional SDE. For nonparametric kernel-based estimators in It\^o's calculus framework, the reader is referred to Kutoyants \cite{KUTOYANTS94} and \cite{KUTOYANTS04}.
\\
\\
The present paper deals with nonparametric estimators of $b$, computed from $N$ independent continuous time observations of the solution of Equation (\ref{main_equation}) on $[0,T]$. Let us mention that it became usual that such functional data are available and can be processed thanks to the improvements of computers. The question of nonparametric drift estimation in stochastic differential equations from such data has been studied in Comte and Genon-Catalot \cite{CGC19} who consider an It\^o's calculus framework.\\
On the one hand, we extend the functional least squares strategy of Comte and Genon-Catalot \cite{CGC19} to fractional SDE by replacing It\^o's integral by Skorokhod's integral in the objective function defining their projection estimator. The Skorokhod integral is a Malliavin calculus based stochastic integral extending It\^o's integral to nonsemimartingale Gaussian signals as the fractional Brownian motion, but which has the major drawback to be non computable in general. For this reason, the oracle $\widehat b$ mentioned above is just an auxiliary {\it entity}, on which we are able to establish a satisfactory risk bound. Thanks to the relationship between the solution $X_{x_0}$ of Equation (\ref{main_equation}) and its derivative with respect to the initial condition $x_0$, we define an approximate estimator $\widehat b_{\varepsilon}$ of $\widehat b$, calculable from $N$ i.i.d. paths of the couple $(X_{x_0},X_{x_0 +\varepsilon})$, where $x_0$ and $x_0 +\varepsilon$ are two close initial conditions. We deduce a risk bound on the estimator $\widehat b_{\varepsilon}$ from the risk bound established on $\widehat b$.\\
On the other hand, by using the relationship between $X_{x_0}$ and $\partial_{x_0}X_{x_0}$ directly, we define an estimator $\widehat b_{\varepsilon}^{\dag,\prime}$ of the derivative $b'$ of $b$, also calculable from $N$ i.i.d. paths of the couple $(X_{x_0},X_{x_0 +\varepsilon})$. We prove a risk bound on $\widehat b_{\varepsilon}^{\dag,\prime}$ with a parametric rate of convergence, but then, we show that deducing a risk bound on a primitive estimator $\widehat b_{\varepsilon}^{\dag}$ of $b$ is not straightforward unless the function is known at one point. In addition, a compactness condition on the support of $b$ must be added. For this reason, the estimator $\widehat b_{\varepsilon}$ is of interest to estimate $b$ and $\widehat b_{\varepsilon}^{\dag, \prime}$ is of interest to estimate $b'$, both functions may be useful depending on the application context.
\\
Note that almost all the references cited above on the statistical inference for fractional SDE are based on long-time behavior properties of the solution which are often difficult to check in practice, but not required here.
\\
\\
The oracle $\widehat b$ and the estimator $\widehat b_{\varepsilon}$ are respectively studied in Subsections \ref{subsection_theoretical_estimator} and \ref{subsection_calculable_estimator} of Section \ref{section_projection_estimators}. Subsection \ref{subsection_rates_usual_bases} provides examples of function bases well adapted in our situation. We can in these frameworks obtain convergence results and rates. The estimators $\widehat b_{\varepsilon}^{\dag,\prime}$ and $\widehat b_{\varepsilon}^{\dag}$ are studied in Section \ref{section_alternative_estimator}. Lastly, concluding remarks are gathered in Section \ref{concluding} while most proofs are postponed in Section \ref{section_proofs}.
\\
\\
\textbf{Notations.} The vector space of Lipschitz continuous maps from $\mathbb R$ into itself is denoted by $\textrm{Lip}(\mathbb R)$ and equipped with the usual Lipschitz semi-norm $\|.\|_{\textrm{Lip}}$. Now, consider $m\in\mathbb N^*$. The Euclidean norm on $\mathbb R^m$ is denoted by $\|.\|_{2,m}$,
\begin{displaymath}
C_{\mathbf b}^{m}(\mathbb R) :=
\left\{\varphi\in C^m(\mathbb R) :
\max_{k\in\llbracket 0,m\rrbracket}
\|\varphi^{(k)}\|_{\infty} <\infty\right\}
\end{displaymath}
and
\begin{displaymath}
\textrm{Lip}_{\mathbf b}^{m}(\mathbb R) :=
\left\{\varphi\in C^m(\mathbb R) :
\|\varphi\|_{\textrm{Lip}_{\mathbf b}^{m}} =
\|\varphi\|_{\textrm{Lip}}\vee
\max_{k\in\llbracket 1,m\rrbracket}\|\varphi^{(k)}\|_{\infty} <\infty\right\}.
\end{displaymath}
Note that $C_{\mathbf b}^{m}(\mathbb R)\subset\textrm{Lip}_{\mathbf b}^{m}(\mathbb R)$. Finally, for every $n\in\mathbb N^*$, the vector space of infinitely continuously differentiable maps $f :\mathbb R^n\rightarrow\mathbb R$ such that $f$ and all its partial derivatives have polynomial growth is denoted by $C_{\mathbf p}^{\infty}(\mathbb R^n,\mathbb R)$.
%


%
\section{Stochastic integrals with respect to the fractional Brownian motion}\label{section_preliminaries}
This section presents two different methods to define a stochastic integral with respect to the fractional Brownian motion. The first one is based on the pathwise properties of the fractional Brownian motion. Another stochastic integral with respect to the fractional Brownian motion is defined via the Malliavin divergence operator. This stochastic integral is called Skorokhod's integral with respect to $B$. If $H = 1/2$, which means that $B$ is a Brownian motion, the Skorokhod integral defined via the divergence operator coincides with It\^o's integral on its domain. This integral is appropriate to define a suitable oracle of the drift function $b$ in Equation (\ref{main_equation}), while the first one is used in Section \ref{subsection_calculable_estimator} to propose a calculable approximation.
%


%
\subsection{The pathwise stochastic integral}\label{subsection_pathwise_integral}
This subsection deals with some definitions and basic properties of the pathwise stochastic integral with respect to the fractional Brownian motion of Hurst index greater than $1/2$.
%


%
\begin{definition}\label{Riemann_sum}
Consider $x$ and $w$ two continuous functions from $[0,T]$ into $\mathbb R$. Consider a dissection $D := (t_0,\dots,t_m)$ of $[s,t]$ with $m\in\mathbb N^*$ and $s,t\in [0,T]$ such that $s < t$. The Riemann sum of $x$ with respect to $w$ on $[s,t]$ for the dissection $D$ is
\begin{displaymath}
J_{x,w,D}(s,t) :=
\sum_{k = 0}^{m - 1}x(t_k)(w(t_{k + 1}) - w(t_k)).
\end{displaymath}
\end{definition}
\noindent
\textbf{Notation.} With the notations of Definition \ref{Riemann_sum}, the mesh of the dissection $D$ is
\begin{displaymath}
\pi(D) :=
\max_{k\in\llbracket 0,m - 1\rrbracket}
|t_{k + 1} - t_k|.
\end{displaymath}
The following theorem ensures the existence and the uniqueness of Young's integral (see Friz and Victoir \cite{FV10}, Theorem 6.8).
%


%
\begin{theorem}\label{Young_integral}
Let $x$ (resp. $w$) be a $\alpha$-H\"older (resp. $\beta$-H\"older) continuous map from $[0,T]$ into $\mathbb R$ with $\alpha,\beta\in ]0,1]$ such that $\alpha +\beta > 1$. There exists a unique continuous map $J_{x,w} : [0,T]\rightarrow\mathbb R$ such that for every $s,t\in [0,T]$ satisfying $s < t$ and any sequence $(D_n)_{n\in\mathbb N}$ of dissections of $[s,t]$ such that $\pi(D_n)\rightarrow 0$ as $n\rightarrow\infty$,
\begin{displaymath}
\lim_{n\rightarrow\infty}
|J_{x,w}(t) - J_{x,w}(s) - J_{x,w,D_n}(s,t)| = 0.
\end{displaymath}
The map $J_{x,w}$ is the Young integral of $x$ with respect to $w$ and $J_{x,w}(t) - J_{x,w}(s)$ is denoted by
\begin{displaymath}
\int_{s}^{t}x(u)dw(u)
\end{displaymath}
for every $s,t\in [0,T]$ such that $s < t$.
\end{theorem}
\noindent
For any $\alpha\in ]1/2,H[$, the paths of $B$ are $\alpha$-H\"older continuous (see Nualart \cite{NUALART06}, Section 5.1). So, for every process $Y = (Y(t))_{t\in [0,T]}$ with $\beta$-H\"older continuous paths from $[0,T]$ into $\mathbb R$ such that $\alpha +\beta > 1$, by Theorem \ref{Young_integral}, it is natural to define the pathwise stochastic integral of $Y$ with respect to $B$ by
\begin{displaymath}
\left(\int_{0}^{t}Y(s)dB(s)\right)(\omega) :=
\int_{0}^{t}Y(\omega,s)dB(\omega,s)
\end{displaymath}
for every $\omega\in\Omega$ and $t\in [0,T]$.
%


%
\subsection{Skorokhod's integral and density of the solution}
This subsection deals with some definitions and results on Malliavin calculus.
\\
\\
Let $\mathbb L^0([0,T],\mathbb R)$ be the space of measurable functions from $[0,T]$ into $\mathbb R$, and consider the reproducing kernel Hilbert space
\begin{displaymath}
\mathcal H :=
\{h\in\mathbb L^0([0,T],\mathbb R) :
\langle h,h\rangle_{\mathcal H} <\infty\}
\end{displaymath}
of $B$, where $\langle .,.\rangle_{\mathcal H}$ is the inner product defined by
\begin{displaymath}
\langle h,\eta\rangle_{\mathcal H} :=
H(2H - 1)
\int_{0}^{T}\int_{0}^{T}
|t - s|^{2H - 2}h(s)\eta(s)dsdt
\end{displaymath}
for every $h,\eta\in\mathbb L^0([0,T],\mathbb R)$. Let $(\mathbf B(h))_{h\in\mathcal H}$ be the isonormal Gaussian process defined by
\begin{displaymath}
\mathbf B(h) :=
\int_{0}^{.}h(s)dB(s),
\end{displaymath}
which is the Wiener integral of $h\in\mathcal H$ with respect to $B$.
%


%
\begin{definition}\label{Malliavin_derivative}
The Malliavin derivative of a smooth functional
\begin{displaymath}
F = f(
\mathbf B(h_1),\dots,
\mathbf B(h_n))
\end{displaymath}
where $n\in\mathbb N^*$, $f\in C_{\mathbf p}^{\infty}(\mathbb R^n,\mathbb R)$ and $h_1,\dots,h_n\in\mathcal H$, is the $\mathcal H$-valued random variable
\begin{displaymath}
\mathbf DF :=
\sum_{k = 1}^{n}
\partial_k f
(
\mathbf B(h_1),\dots,
\mathbf B(h_n))h_k.
\end{displaymath}
\end{definition}
\noindent
The key property of the operator $\mathbf D$ is the following one.
%


%
\begin{proposition}\label{Malliavin_derivative_domain}
The map $\mathbf D$ is closable from $\mathbb L^2(\Omega,\mathcal A,\mathbb P)$ into $\mathbb L^2(\Omega;\mathcal H)$. Its domain in $\mathbb L^2(\Omega,\mathcal A,\mathbb P)$, denoted by $\mathbb D^{1,2}$, is the closure of the smooth functionals space for the seminorm $\|.\|_{1,2}$ defined by
\begin{displaymath}
\|F\|_{1,2}^{2} :=
\mathbb E(|F|^2) +
\mathbb E(\|\mathbf DF\|_{\mathcal H}^{2}) < \infty
\end{displaymath}
for every $F\in\mathbb L^2(\Omega,\mathcal A,\mathbb P)$. The Malliavin derivative of $F\in\mathbb D^{1,2}$ at time $s\in [0,T]$ is denoted by $\mathbf D_sF$.
\end{proposition}
\noindent
For a proof, see Nualart \cite{NUALART06}, Proposition 1.2.1.
%


%
\begin{definition}\label{divergence_operator}
The adjoint $\delta$ of the Malliavin derivative $\mathbf D$ is the divergence operator. The domain of $\delta$ is denoted by $\normalfont{\textrm{dom}}(\delta)$, and $u\in\normalfont{\textrm{dom}}(\delta)$ if and only if there exists a deterministic constant $\mathfrak c_u > 0$ such that for every $F\in\mathbb D^{1,2}$,
\begin{displaymath}
|\mathbb E(\langle\mathbf DF,u\rangle_{\mathcal H})|
\leqslant
\mathfrak c_u\mathbb E(|F|^2)^{1/2}.
\end{displaymath}
\end{definition}
\noindent
For any process $Y = (Y(s))_{s\in [0,T]}$ and every $t\in ]0,T]$, if $Y\mathbf 1_{[0,t]}\in\textrm{dom}(\delta)$, then its Skorokhod integral with respect to $B$ is defined on $[0,t]$ by
\begin{displaymath}
\int_{0}^{t}Y(s)\delta B(s) :=
\delta(Y\mathbf 1_{[0,t]}),
\end{displaymath}
and its Skorokhod integral with respect to $X$ is defined by
\begin{displaymath}
\int_{0}^{t}
Y(s)\delta X(s) :=
\int_{0}^{t}Y(s)b(X(s))ds +
\sigma\int_{0}^{t}Y(s)\delta B(s).
\end{displaymath}
Note that since $\delta$ is the adjoint of the Malliavin derivative $\mathbf D$, the Skorokhod integral of $Y$ with respect to $B$ on $[0,t]$ is a centered random variable. Indeed,
\begin{equation}\label{divergence_operator_1}
\mathbb E\left(\int_{0}^{t}Y(s)\delta B(s)\right) =
\mathbb E(1\cdot\delta(Y\mathbf 1_{[0,t]})) =
\mathbb E(\langle\mathbf D(1),Y\mathbf 1_{[0,t]}\rangle_{\mathcal H}) =
0.
\end{equation}
Let $\mathcal S$ be the space of the smooth functionals presented in Definition \ref{Malliavin_derivative} and consider $\mathbb D^{1,2}(\mathcal H)$, the closure of
\begin{displaymath}
\mathcal S_{\mathcal H} :=
\left\{
\sum_{j = 1}^{n}F_jh_j
\textrm{ $;$ }
h_1,\dots,h_n\in\mathcal H
\textrm{, }
F_1,\dots,F_n\in\mathcal S
\right\}
\end{displaymath}
for the seminorm $\|.\|_{1,2,\mathcal H}$ defined by
\begin{displaymath}
\|u\|_{1,2,\mathcal H}^{2} :=\mathbb E(\|u\|_{\mathcal H}^{2}) +
\mathbb E(\|\mathbf Du\|_{\mathcal H\otimes\mathcal H}^{2}) <\infty
\end{displaymath}
for every $u\in\mathbb L^2(\Omega\times [0,T])$ (see Nualart \cite{NUALART06}, p. 31, Remark 2). The following proposition provides an isometry type property for the Skorokhod integral with respect to $B$ on $\mathbb D^{1,2}(\mathcal H)$, which is a subspace of $\textrm{dom}(\delta)$ by Nualart \cite{NUALART06}, Proposition 1.3.1. This result is useful for our purpose and is proved in Biagini et al. \cite{BHOZ08} (see Theorem 3.11.1).
%


%
\begin{proposition}\label{isometry_Skorokhod}
For every $Y,Z\in\mathbb D^{1,2}(\mathcal H)$,
\begin{eqnarray*}
 \mathbb E(\delta(Y)\delta(Z)) & = &
 \alpha_H\int_{0}^{T}\int_{0}^{T}
 \mathbb E(Y(u)Z(v))|v - u|^{2H - 2}dvdu\\
 & & +
 \alpha_{H}^{2}\int_{[0,T]^4}
 \mathbb E(\mathbf D_{u'}Y(v)\mathbf D_{v'}Z(u))|u - u'|^{2H - 2}|v - v'|^{2H - 2}dudu'dvdv'.
\end{eqnarray*}
\end{proposition}
\noindent
In the sequel, the function $b$ fulfills the following assumption.
%


%
\begin{assumption}\label{assumption_b}
The function $b$ belongs to $C^1(\mathbb R)$ and there exist $m,M\in\mathbb R$ such that
\begin{displaymath}
m\leqslant b'(x)\leqslant M
\textrm{ $;$ }
\forall x\in\mathbb R.
\end{displaymath}
\end{assumption}
\noindent
Under Assumption \ref{assumption_b}, the following result is a straightforward application of Proposition \ref{isometry_Skorokhod} to functionals of the solution $X$ of Equation (\ref{main_equation}).
%


%
\begin{corollary}\label{isometry_Skorokhod_fX}
Let $X$ be the solution of Equation (\ref{main_equation}).
Under Assumption \ref{assumption_b}, $X\in\mathbb D^{1,2}(\mathcal H)$ and for every $\varphi,\psi\in\normalfont{\textrm{Lip}}_{\mathbf b}^{1}(\mathbb R)$,
\begin{displaymath}
\mathbb E(\delta(\varphi(X))\delta(\psi(X))) =
\alpha_H\int_{0}^{T}\int_{0}^{T}
\mathbb E(\varphi(X(u))\psi(X(v)))|v - u|^{2H - 2}dvdu + R_{\varphi,\psi}
\end{displaymath}
where
\begin{eqnarray*}
 R_{\varphi,\psi} & := &
 \alpha_{H}^{2}
 \int_{[0,T]^4}
 \mathbb E(\varphi'(X(v))\psi'(X(v'))
 \mathbf D_{u'}X(v)\mathbf D_uX(v'))
 |u - u'|^{2H - 2}|v - v'|^{2H - 2}
 dudu'dvdv'\\
 & = &
 \alpha_{H}^{2}\sigma^2
 \int_{[0,T]^2}
 \int_{0}^{v}\int_{0}^{v'}|u - u'|^{2H - 2}|v - v'|^{2H - 2}\\
 & &
 \times
 \mathbb E\left(\varphi'(X(v))\psi'(X(v'))\exp\left(
 \int_{u'}^{v}b'(X(s))ds +
 \int_{u}^{v'}b'(X(s))ds\right)\right)dudu'dvdv'.
\end{eqnarray*}
\end{corollary}
\noindent
The following theorem provides suitable controls of the moments of Skorokhod's integral.
%


%
\begin{theorem}\label{control_divergence_integral}
Under Assumption \ref{assumption_b}, for every $p > 1/H$, there exists a deterministic constant $\mathfrak c_{p,H,\sigma} > 0$, only depending on $p$, $H$ and $\sigma$, such that for every $\varphi\in\normalfont{\textrm{Lip}}_{\mathbf b}^{1}(\mathbb R)$,
\begin{eqnarray*}
 \mathbb E\left(
 \left|\int_{0}^{T}\varphi(X(s))\delta B(s)\right|^p\right)
 & \leqslant &
 \mathfrak c_{p,H,\sigma}\overline{\mathfrak m}_{p,H,M}(T)
 \left[\left(\int_{0}^{T}\mathbb E(|\varphi(X(s))|^{1/H})ds\right)^{pH}\right.\\
 & &
 +\left.
 \left(\int_{0}^{T}\mathbb E(|\varphi'(X(s))|^p)^{1/(pH)}ds\right)^{pH}\right] <\infty
\end{eqnarray*}
where $\overline{\mathfrak m}_{p,H,M}(T) :=\mathfrak m_{p,H,M}(T)\vee 1$ and
\begin{displaymath}
\mathfrak m_{p,H,M}(T) :=
\left(-\frac{H}{M}\right)^{pH}\mathbf 1_{M < 0} +
T^{pH}\mathbf 1_{M = 0} +
\left(\frac{H}{M}\right)^{pH}e^{pMT}\mathbf 1_{M > 0}.
\end{displaymath}
\end{theorem}
\noindent
Note that if $M < 0$, then Theorem \ref{control_divergence_integral} has been already proved in Hu, Nualart and Zhou \cite{HNZ18} (see Proposition 4.4.(2)).
%


%
\begin{remark}\label{variance_Skorokhod_integral}
On the one hand, note that the control of the variance of Skorokhod's integral provided in Theorem \ref{control_divergence_integral} is a straightforward consequence of Corollary \ref{isometry_Skorokhod_fX}. On the other hand, with the notations of Corollary \ref{isometry_Skorokhod_fX}, note that for $H = 1/2$, the solution $X$ of Equation (\ref{main_equation}) is adapted and then
\begin{displaymath}
R_{\varphi,\psi} =
\int_{0}^{T}\int_{0}^{T}\mathbb E(\mathbf D_uX(v)\mathbf D_vX(u))dudv = 0.
\end{displaymath}
This reduces importantly the order of the variance of Skorokhod's integral with respect to the case $H > 1/2$.
\end{remark}
\noindent
Lastly, the following proposition provides an expression and a bound for the density of the solution to Equation (\ref{main_equation}).
%


%
\begin{proposition}\label{density_solution_main_equation}
Under Assumption \ref{assumption_b}, for any $t\in ]0,T]$, the probability distribution of $X^*(t) := X(t) -\mathbb E(X(t))$ has a $\mathbb R$-supported density with respect to Lebesgue's measure $p_{t}^{*}(x_0,.)$ such that, for every $x\in\mathbb R$,
\begin{displaymath}
p_{t}^{*}(x_0,x) =
\frac{\mathbb E(|X^*(t)|)}{2g_{t}^{*}(x_0,x)}
\exp\left(-\int_{0}^{x}\frac{zdz}{g_{t}^{*}(x_0,z)}\right)
\end{displaymath}
where
\begin{eqnarray*}
 g_{t}^{*}(x_0,x) & := &
 \mathbb E(\langle\mathbf DX(t),-\mathbf D\mathbf L^{-1}X(t)\rangle_{\mathcal H}|X^*(t) = x)\\
 & &
 \hspace{4cm}\in
 [\sigma(m,t),\sigma(M,t)]\subset ]0,\infty[,
\end{eqnarray*}
$\mathbf L$ is the Ornstein-Uhlenbeck operator and
\begin{displaymath}
\sigma(\mu,t)^2 :=
\alpha_H\sigma^2\int_{0}^{t}\int_{0}^{t}|v - u|^{2H - 2}
e^{\mu(2t - v - u)}dudv > 0
\textrm{ $;$ }
\forall\mu\in\mathbb R.
\end{displaymath}
In particular, for every $x\in\mathbb R$,
\begin{displaymath}
\frac{\mathbb E(|X^*(t)|)}{2\sigma(M,t)^2}\exp\left(-\frac{x^2}{2\sigma(m,t)^2}\right)
\leqslant p_{t}^{*}(x_0,x)\leqslant
\frac{\mathbb E(|X^*(t)|)}{2\sigma(m,t)^2}\exp\left(-\frac{x^2}{2\sigma(M,t)^2}\right).
\end{displaymath}
\end{proposition}
\noindent
A straightforward consequence of Proposition \ref{density_solution_main_equation} is that for any $t\in ]0,T]$, the probability distribution of $X(t)$ has a $\mathbb R$-supported density with respect to Lebesgue's measure $p_t(x_0,.)$ such that, for every $x\in\mathbb R$,
\begin{displaymath}
p_t(x_0,x) =
p_{t}^{*}(x_0,x -\mathbb E(|X(t)|))
\end{displaymath}
and
\begin{equation}\label{control_density_solution}
\mathfrak m_t(x_0,x,M,m)
\leqslant p_t(x_0,x)\leqslant
\mathfrak m_t(x_0,x,m,M),
\end{equation}
where
\begin{displaymath}
\mathfrak m_t(x_0,x,\mu_1,\mu_2) :=
\frac{\mathbb E(|X^*(t)|)}{2\sigma(\mu_1,t)^2}\exp\left[-\frac{(x -\mathbb E(|X(t)|))^2}{2\sigma(\mu_2,t)^2}\right]
\textrm{ $;$ }
\forall\mu_1,\mu_2\in\mathbb R.
\end{displaymath}
Since the paths of $X$ are $\alpha$-H\"older continuous for any $\alpha\in ]0,H[$,
\begin{eqnarray*}
 \mathbb E(|X^*(t)|) & = &
 \mathbb E(|X(t) - x_0 -\mathbb E(X(t) - x_0)|)\\
 & \leqslant &
 2\mathbb E(|X(t) - X(0)|)
 \leqslant
 2\mathbb E(\|X\|_{\alpha\textrm{-H\"ol},T})t^{\alpha}
\end{eqnarray*}
where
\begin{displaymath}
\|X\|_{\alpha\textrm{-H\"ol},T} :=
\sup_{0\leqslant s < t\leqslant T}
\frac{|X(t) - X(s)|}{|t - s|^{\alpha}},
\end{displaymath}
which has a finite first order moment because $\mathbb E(\|B\|_{\alpha\textrm{-H\"ol},T}) <\infty$ and $b$ is Lipschitz continuous. Then, since $\sigma(m,t)^2\geqslant\sigma^2e^{-2\|b'\|_{\infty}T}t^{2H}$,
\begin{displaymath}
\mathfrak m_t(x_0,x,m,M)\leqslant
\mathfrak c_Tt^{\alpha -2H}
\end{displaymath}
with
\begin{displaymath}
\mathfrak c_T :=
\frac{\mathbb E(\|X\|_{\alpha\textrm{-H\"ol},T})}{\sigma^2e^{-2\|b'\|_{\infty}T}}.
\end{displaymath}
Therefore, by taking $\alpha\in ]2H - 1,H[$, Inequality (\ref{control_density_solution}) implies that for every $x\in\mathbb R$, $p_.(x_0,x)\in\mathbb L^1(]0,T],dt)$.
%


%
\section{Projection estimator of the drift function}\label{section_projection_estimators}
Under Assumption \ref{assumption_b}, $b$ is Lipschitz continuous on $\mathbb R$ and its derivative is bounded. So, Equation (\ref{main_equation}) has a unique solution $X$ and the associated It\^o map $\mathcal I$ is continuously differentiable from $\mathbb R\times C^0([0,T],\mathbb R)$ into $C^0([0,T],\mathbb R)$.
%


%
\subsection{A Skorokhod's integral based oracle}\label{subsection_theoretical_estimator}
This subsection deals with an oracle of $b$, constructed as the estimator of Comte and Genon-Catalot \cite{CGC19} by replacing It\^o's integral by Skorokhod's one. The risk bound on this oracle will allow us to establish a risk bound on an estimator at Subsection \ref{subsection_calculable_estimator}.
%


%
\subsubsection{The objective function}\label{subsubsection_objective_function}
Let $f_T$ be the function defined by
\begin{displaymath}
f_T(x) :=\frac{1}{T}\int_{0}^{T}p_s(x_0,x)ds
\textrm{ $;$ }
\forall x\in\mathbb R,
\end{displaymath}
where $p_s(x_0,.)$ is the density with respect to Lebesgue's measure of the probability distribution of $X(s)$ for any $s\in ]0,T]$. For any $x\in\mathbb R$, $p_s(x_0,x)$ is well defined and belongs to $]0,\infty[$ because ${\rm supp}(p_s(x_0,.)) =\mathbb R$ (see Proposition \ref{density_solution_main_equation}). So, since $p_.(x_0,x)\in\mathbb L^1(]0,T],dt)$ as established at the end of Section \ref{section_preliminaries}, $f_T(x)$ is well defined, belongs to $]0,\infty[$, and then ${\rm supp}(f_T) =\mathbb R$. Moreover, by Fubini-Tonelli's theorem, $f_T$ is measurable and
\begin{displaymath}
\int_{-\infty}^{\infty}f_T(x)dx =
\frac{1}{T}\int_{0}^{T}\int_{-\infty}^{\infty}p_s(x_0,x)dxds = 1.
\end{displaymath}
So, $f_T$ is a $\mathbb R$-supported density function.
\\
\\
Now, consider $N\in\mathbb N^*$ independent copies $B_1,\dots,B_N$ of $B$, $X^i :=\mathcal I(x_0,B^i)$ for every $i\in\{1,\dots,N\}$, and the objective function $\gamma_N$ defined by
\begin{displaymath}
\gamma_N(\tau) :=
\frac{1}{NT}\sum_{i = 1}^{N}\left(
\int_{0}^{T}\tau(X^i(s))^2ds - 2\int_{0}^{T}\tau(X^i(s))\delta X^i(s)\right)
\end{displaymath}
for every function $\tau :\mathbb R\rightarrow\mathbb R$.
\\
\\
Note that for any bounded function $\tau$ from $\mathbb R$ into itself, thanks to Equality (\ref{divergence_operator_1}),
\begin{eqnarray*}
 \mathbb E(\gamma_N(\tau)) & = &
 \frac{1}{T}\int_{0}^{T}\mathbb E(\tau(X(s))^2 -2\tau(X(s))b(X(s)))ds +
 \frac{\sigma}{T}\mathbb E\left(\int_{0}^{T}\tau(X(s))\delta B(s)\right)\\
 & = &
 \frac{1}{T}\int_{0}^{T}\mathbb E((\tau(X(s)) - b(X(s)))^2)ds -\frac{1}{T}\int_{0}^{T}\mathbb E(b(X(s))^2)ds.
\end{eqnarray*}
Then, the definition of $f_T$ gives
\begin{equation}\label{EspGamma}
\mathbb E(\gamma_N(\tau)) =
\int_{-\infty}^{\infty}(\tau(x) - b(x))^2f_T(x)dx -
\int_{-\infty}^{\infty}b(x)^2f_T(x)dx.
\end{equation}
Equality (\ref{EspGamma}) shows that $\mathbb E(\gamma_N(\tau))$ is the smallest for $\tau$ the nearest of $b$. Therefore, to minimize its empirical version $\gamma_N(\tau)$ should provide a function near of $b$.
%


%
\begin{remark}\label{why_not_pathwise}
The pathwise stochastic integral with respect to $B$, defined in Subsection \ref{subsection_pathwise_integral}, is not centered in general, and not even  for $H=1/2$. Indeed, if $H = 1/2$, then it coincides with Stratonovich's integral. This is the main reason why the objective function above is defined via Skorokhod's integral. Moreover, the pathwise stochastic integral doesn't satisfy an isometry type property as Skorokhod's integral (see Proposition \ref{isometry_Skorokhod}), which is crucial in the sequel.
\end{remark}
%


%
\subsubsection{The oracle}\label{subsubsection_projection_estimators}
Consider $A\in\mathcal B(\mathbb R)$ and assume that $\mathbb L^2(A,dx)$ (resp. $\mathbb L^2(A,f_T(x)dx)$) is equipped with its usual inner product $\langle .,.\rangle$ (resp. $\langle .,.\rangle_{f_T}$). For any $m\in\mathbb N^*$, consider also
\begin{displaymath}
\mathcal S_m :=
\textrm{span}\{\varphi_0,\dots,\varphi_{m - 1}\},
\end{displaymath}
where $(\varphi_0,\dots,\varphi_{m - 1})$ is an orthonormal family of $\mathbb L^2(A,dx)$. Moreover, assume that the functions $\varphi_j$, $j\in\mathbb N$ are bounded. So, $\mathcal S_m\subset\mathbb L^2(A,f_T(x)dx)$.
\\
\\
Consider
\begin{displaymath}
\widehat{\mathbf\Psi}(m) :=
\left(\frac{1}{NT}
\sum_{i = 1}^{N}
\int_{0}^{T}\varphi_j(X^i(s))\varphi_k(X^i(s))ds\right)_{j,k = 0,\dots,m - 1},
\end{displaymath}
assume that this matrix is invertible, and let
\begin{displaymath}
\widehat b_m =\arg\min_{\tau\in\mathcal S_m}\gamma_N(\tau)
\end{displaymath}
be the Skorokhod's integral based projection oracle of $b_A := b_{|A}$ on $\mathcal S_m$. As in Comte and Genon-Catalot \cite{CGC19}, Section 2.2,
\begin{displaymath}
\widehat b_m =
\sum_{j = 0}^{m - 1}\widehat\theta_j\varphi_j
\end{displaymath}
where
\begin{displaymath}
\widehat{\theta}(m) :=
(\widehat\theta_0,\dots,\widehat\theta_{m - 1})^* =
\widehat{\mathbf\Psi}(m)^{-1}\widehat{\mathbf x}(m)
\end{displaymath}
with
\begin{displaymath}
\widehat{\mathbf x}(m) :=
\left(\frac{1}{NT}
\sum_{i = 1}^{N}
\int_{0}^{T}\varphi_j(X^i(s))\delta X^i(s)\right)_{j = 0,\dots,m - 1}^{*}.
\end{displaymath}
Note that
\begin{displaymath}
\widehat{\mathbf x}(m) =
(\langle\varphi_j,b\rangle_N)_{j = 0,\dots,m - 1}^{*} +
\mathbf e(m)
\end{displaymath}
and
\begin{displaymath}
\widehat{\mathbf\Psi}(m) =
(\langle\varphi_j,\varphi_k\rangle_N)_{j,k = 0,\dots,m - 1},
\end{displaymath}
where
\begin{displaymath}
\langle\tau_1,\tau_2\rangle_N :=
\frac{1}{NT}\sum_{i = 1}^{N}\int_{0}^{T}\tau_1(X^i(s))\tau_2(X^i(s))ds
\end{displaymath}
for every measurable functions $\tau_1,\tau_2 :\mathbb R\rightarrow\mathbb R$ (the associated norm is denoted by $\|.\|_N$), and
\begin{displaymath}
\mathbf e(m) :=
\left(\frac{\sigma}{NT}\sum_{i = 1}^{N}\int_{0}^{T}\varphi_j(X^i(s))\delta B^i(s)\right)_{j = 0,\dots,m - 1}^{*}.
\end{displaymath}
By Equality (\ref{divergence_operator_1}), $\mathbf e(m)$ is centered, as expected for an error term in regression.
%


%
\subsubsection{Risk bounds on the oracle}\label{subsubsection_convergence_results}
Throughout this subsection, $f_T$ and the functions $\varphi_j$, $j\in\mathbb N$ fulfill the following assumption.
%


%
\begin{assumption}\label{assumption_basis_1}
$\lambda(A) > 0$ and, for $m\leqslant NT$,
\begin{enumerate}
 \item $(\varphi_0,\dots,\varphi_{m - 1})$ is an orthonormal family of $\mathbb L^2(A,dx)$.
 \item The functions $\varphi_j$, $j = 0,\dots,m - 1$, are bounded and belong to $C_{\mathbf b}^{1}(A)$.
 \item There exist $x_0,\dots,x_{m - 1}\in A$ such that
 \begin{displaymath}
 \det[(\varphi_j(x_k))_{j,k = 0,\dots,m - 1}]\not= 0.
 \end{displaymath}
\end{enumerate}
\end{assumption} 
\noindent
By Comte and Genon-Catalot \cite{CGC19}, Lemma 1, which remains true for $H > 1/2$ without additional arguments,
\begin{displaymath}
\mathbf\Psi(m) :=\mathbb E(\widehat{\mathbf\Psi}(m))= \left(\int_A\varphi_j(x)\varphi_k(x)f_T(x)dx\right)_{0\leqslant j, k \leqslant m-1}
\end{displaymath}
is invertible under Assumption \ref{assumption_basis_1}. In addition, we impose that 
\begin{displaymath}
L(m) :=
\sup_{x\in A}
\sum_{j = 0}^{m - 1}
\varphi_j(x)^2
\textrm{ and }
 R(m) :=
\sup_{x\in A}\sum_{j = 0}^{m - 1}
\varphi_j'(x)^2
\end{displaymath}
fulfill the following assumption.
%


%
\begin{assumption}\label{assumption_basis_2}
There exists $\rho > 0$ and $\kappa\geqslant 1$ such that $R(m)\leqslant\rho L(m)^{\kappa}$ and
\begin{displaymath}
L(m)
(\|\mathbf\Psi(m)^{-1}\|_{\normalfont\textrm{op}}\vee 1)
\leqslant
\frac{\mathfrak c_{\kappa,T}}{2}\cdot\frac{NT}{\log(NT)}
\quad
\textrm{with}
\quad
\mathfrak c_{\kappa,T} :=\frac{3\log(3/2)-1}{(7 +\kappa)T}.
\end{displaymath}
\end{assumption}
\noindent
The above condition is a generalization of the so-called {\it stability condition} introduced for the standard regression by Cohen et al. \cite{CDL13,CDL19}, also considered in Comte and Genon-Catalot \cite{CGC19}.
\\
\\
\textbf{Convention.} When ${\bf M}$ is a symmetric nonnegative and noninvertible matrix, $\|{\bf M}^{-1}\|_{\rm op} :=\infty$. This is a coherent convention because if ${\bf M}$ is invertible, then $\|{\bf M}^{-1}\|_{\rm op} = 1/\inf\{{\rm sp}({\bf M})\}$.
\\
\\
In order to ensure the existence of the oracle and ti be able to bound its integrated risk, $\widehat b_m$ is replaced by
\begin{displaymath}
\widetilde b_m :=
\widehat b_m\mathbf 1_{\widehat\Lambda_{\kappa}(m)},
\end{displaymath}
where
\begin{equation}\label{Lambdakappa}
\widehat\Lambda_{\kappa}(m) :=
\left\{L(m)(\|\widehat{\mathbf\Psi}(m)^{-1}\|_{\textrm{op}}\vee 1)
\leqslant\mathfrak c_{\kappa,T}\frac{NT}{\log(NT)}\right\}.
\end{equation}
Note that with the previous convention, on the event $\widehat\Lambda_{\kappa}(m)$, $\widehat{\mathbf\Psi}(m)$ is invertible because
\begin{displaymath}
\inf\{{\rm sp}(\widehat{\mathbf\Psi}(m))\}
\geqslant
\frac{L(m)}{\mathfrak c_{\kappa,T}}\cdot
\frac{\log(NT)}{NT}.
\end{displaymath}
Then, $\widetilde b_m$ is well-defined. Moreover, necessarily, $m\leqslant NT/\log(NT)$ on $\widehat\Lambda_{\kappa}(m)$. 
\\
\\
The two following results provide controls of the empirical risk and of the $f_T$-weighted integrated  risk of $\widetilde b_m$ respectively.
%


%
\begin{theorem}\label{empirical_risk_bound}
Under Assumptions \ref{assumption_b}, \ref{assumption_basis_1} and \ref{assumption_basis_2},
\begin{displaymath}
\mathbb E(\|\widetilde b_m - b_A\|_{N}^{2})
\leqslant
\inf_{\tau\in\mathcal S_m}
\|\tau - b_A\|_{f_T}^{2} +
\frac{2}{NT}\normalfont{\textrm{trace}}(\mathbf\Psi(m)^{-1}\mathbf\Psi(m,\sigma)) +
\mathfrak c_{\rho,\kappa,\sigma}(1 +\mathfrak b_T)\frac{\overline{\mathfrak m}_{2,H,M}(T)}{NT}
\end{displaymath}
where $\mathfrak c_{\rho,\kappa,\sigma} > 0$ is a deterministic constant depending only on $\rho$, $\kappa$ and $\sigma$, $\overline{\mathfrak m}_{2,H,M}(T)$ is a constant defined in Theorem \ref{control_divergence_integral},
\begin{displaymath}
\mathbf\Psi(m,\sigma) :=
\left(\frac{\sigma^2}{T}
\mathbb E\left(
\left(\int_{0}^{T}\varphi_j(X(s))\delta B(s)\right)
\left(\int_{0}^{T}\varphi_k(X(s))\delta B(s)\right)\right)\right)_{j,k = 0,\dots,m - 1}
\end{displaymath}
and $\mathfrak b_T :=\|b_{A}^{2}\|_{f_T}$.
\end{theorem}
%


%
\begin{corollary}\label{risk_bound}
Under Assumptions \ref{assumption_b}, \ref{assumption_basis_1} and \ref{assumption_basis_2},
\begin{eqnarray*}
 \mathbb E(\|\widetilde b_m - b_A\|_{f_T}^{2})
 & \leqslant &
 \left(1 +  \frac{4T\mathfrak c_{\kappa,T}}{\log(NT)}\right)
 \inf_{\tau\in\mathcal S_m}
 \|\tau - b_A\|_{f_T}^{2}\\
 & &
 \hspace{1cm}
 +\frac{8}{NT}\normalfont{\textrm{trace}}(\mathbf\Psi(m)^{-1}\mathbf\Psi(m,\sigma))
 +\overline{\mathfrak c}_{\rho,\kappa,\sigma}(1 +\mathfrak b_T)
 \frac{\overline{\mathfrak m}_{2,H,M}(T)}{NT}
\end{eqnarray*}
where $\overline{\mathfrak c}_{\rho,\kappa,\sigma} > 0$ is a deterministic constant depending only on $\rho$, $\kappa$ and $\sigma$.
\end{corollary}
%


%
\begin{remark}\label{remark_variance_term}
Note that
\begin{eqnarray*}
 \normalfont{\textrm{trace}}(\mathbf\Psi(m)^{-1}\mathbf\Psi(m,\sigma))
 & = &
 \normalfont{\textrm{trace}}(\mathbf\Psi(m)^{-1/2}\mathbf\Psi(m,\sigma)\mathbf\Psi(m)^{-1/2})
 \leqslant
 m\|\mathbf\Psi(m)^{-1/2}\mathbf\Psi(m,\sigma)\mathbf\Psi(m)^{-1/2}\|_{\normalfont{\textrm{op}}}\\
 & = &
 m\sup_{\tau\in S_m :\|\tau\|_{f_T} = 1}
 \mathbb E\left[\left(\int_{0}^{T}\tau(X(s))\delta B(s)\right)^2\right].
\end{eqnarray*}
\end{remark}
\noindent
The risk decompositions given in Theorem \ref{empirical_risk_bound} and Corollary \ref{risk_bound} both involve the same types of terms:
\begin{itemize}
 \item The first one is equal or proportional to $\inf_{\tau\in\mathcal S_m}\|\tau - b_A\|_{f_T}^{2}$ and is a squared bias term due to the projection strategy. It is decreasing when $m$ increases, because then the projection space grows.
 \item The second one, $\normalfont{\textrm{trace}}(\mathbf\Psi(m)^{-1}\mathbf\Psi(m,\sigma))/(NT)$, is a variance term. From the remark above, it is bounded by $m\|\mathbf\Psi(m)^{-1}\mathbf\Psi(m,\sigma)\|_{\textrm{op}}/(NT)$ which is increasing with $m$.
 \item The last one is a residual negligible term, which is small when  $N$ is large. Note that if the upper-bound $M$ on $b'$ is nonnegative, then $\overline{\mathfrak m}_{2,H,M}(T)$ explodes for large values of $T$.
\end{itemize}
The  order of the bias generally depends on the regularity of the function, and the order of the trace term is discussed below. Both quantities imply that a choice of $m$ ensuring a  compromise between the bias and the variance is required to obtain the convergence of $\widetilde b_m$ and a rate.
\\
\\
Finally, let us provide a control for $\textrm{trace}(\mathbf\Psi(m)^{-1}\mathbf\Psi(m,\sigma))$ which allows comparison with non fractional results.
%


%
\begin{proposition}\label{bound_trace_term}
Under Assumptions \ref{assumption_b} and \ref{assumption_basis_1},
\begin{eqnarray*}
 \frac{
 \normalfont{\textrm{trace}}(\mathbf\Psi(m)^{-1}\mathbf\Psi(m,\sigma))}{NT}
 & \leqslant &
 \mathfrak c_{2,H,\sigma}\sigma^2\frac{\overline{\mathfrak m}_{2,H,M}(T)}{NT^{2 - 2H}}\\
 & &
 \hspace{1.5cm}
 \times\min\{
 (L(m) + R(m))\|\mathbf\Psi(m)^{-1}\|_{{\rm op}}, m(1+R(m)\|\mathbf\Psi(m)^{-1}\|_{{\rm op}})\}.
\end{eqnarray*}
\end{proposition}
\noindent
In the standard case, with $H = 1/2$ and a constant volatility function $\sigma$, it holds that
\begin{displaymath}
\frac{1}{NT}\textrm{trace}(\mathbf\Psi(m)^{-1}\mathbf\Psi(m,\sigma)) =
\sigma^2\frac{m}{NT}
\end{displaymath}
as established in Comte and Genon-Catalot \cite{CGC19}. Here, for $M < 0$, the constant $\overline{\mathfrak m}_{2,H,M}(T)$ does not depend on $T$ and thus, $NT$ becomes $NT^{2-2H}$ which is coherent. However, the additional term $R(m)\|\mathbf\Psi(m)^{-1}\|_{\textrm{op}}$ may have an important order in $m$ and substantially increases the variance. Thus, it will deteriorate the rate of the estimators. So, there is a \textit{discontinuity} between the cases $H = 1/2$ and $H > 1/2$, which is explained in Remark \ref{variance_Skorokhod_integral}. However, note that this discontinuity is specific to the estimation strategy investigated in this paper.
%


%
\subsection{Rates in some usual bases}\label{subsection_rates_usual_bases}
 Now, for projection estimators, different bases can be considered. In the present setting, the bases have to be differentiable. Let us present two examples.
%


%
\subsubsection{Rates on Fourier-Sobolev spaces for the trigonometric basis}\label{subsubsection_discussion_Fourier}
A first example  is the compactly supported trigonometric basis. For $A = [\ell,\texttt r]$, it is defined by
\begin{eqnarray*}
\varphi_0(x) & := & \frac{1}{\sqrt{\texttt r - \ell}}
\mathbf 1_{[\ell,\texttt r]}(x)
\textrm{, }\\
\varphi_{2j+1}(x) & := & \sqrt{\frac{2}{\texttt r- \ell}}\cos\left(2\pi j\frac{x - \ell}{\texttt r - \ell}\right)
\mathbf 1_{[\ell, \texttt r]}(x)\textrm{ and}\\
\varphi_{2j}(x) & := & \sqrt{\frac{2}{\texttt r -\ell}}
\sin\left(2\pi j\frac{x - \ell}{\texttt r - \ell}\right)
\mathbf 1_{[a,b]}(x)
\end{eqnarray*}
for every $x\in\mathbb R$ and $j\geqslant 1$. This basis satisfies, for $m$ odd and any $x\in [\ell,\texttt r]$,
\begin{displaymath}
\sum_{j = 0}^{m - 1}\varphi_{j}^{2}(x) = m
\textrm{ and }
\sup_{x\in [\ell, \texttt r]}
\sum_{j = 0}^{m - 1}
\varphi_j'(x)^2\leqslant
\frac{(2\pi)^2}{(\texttt r - \ell)^3}m^3.
\end{displaymath}
So,
\begin{displaymath}
L(m) = m
\textrm{ and }
R(m) = \rho(\ell,\texttt r)m^3
\end{displaymath}
where $\rho(\ell, \texttt r) = (2\pi)^2/(\texttt r- \ell)^3$.
\\
\\
In the Brownian setting, where $H = 1/2$, for a constant volatility function $\sigma(x)\equiv\sigma$, as recalled above, the variance term is $\sigma^2m/(NT)$ (see Comte and Genon-Catalot \cite{CGC19}). Here, if we assume that $f_T$ is lower bounded on $A$ by $f_0 > 0$, then $\|\mathbf\Psi(m)^{-1}\|_{\textrm{op}}\leqslant 1/f_0$ and the bound of Proposition \ref{bound_trace_term} becomes
\begin{displaymath}
\frac{1}{NT}
\textrm{trace}(\mathbf\Psi(m)^{-1}\mathbf\Psi(m,\sigma))
\leqslant
\mathfrak c_{2,H,\sigma}\sigma^2
\frac{\overline{\mathfrak m}_{2,H,M}(T)}{NT^{2 - 2H}}\cdot\frac{m}{f_0}[1 +\rho(\ell, \texttt r)m^2].
\end{displaymath}
The additional term $R(m) \|\mathbf\Psi(m)^{-1}\|_{{\rm op}}$ discussed after Proposition \ref{bound_trace_term} is here of order $m^3$.
\\
\\
Now, let us evaluate the bias term. Consider $\beta\in\mathbb N^*$ and the Sobolev space
\begin{displaymath}
W_{2}^{\beta}([\ell,\texttt r]) :=
\left\{\varphi : [\ell, \texttt r]\rightarrow\mathbb R :
\int_{\ell}^{\texttt r}|\varphi^{(\beta)}(x)|^2dx <\infty\right\}.
\end{displaymath}
If $b_A\in W_{2}^{\beta}([\ell,\texttt r])$, by DeVore and Lorentz \cite{DL93}, Theorem 2.3 p. 205, then there exists a deterministic constant $\mathfrak c_{\beta,\ell,\texttt r} > 0$, not depending on $m$, such that
\begin{displaymath}
\|p_{\mathcal S_m}^{\perp}(b_A) - b_A\|^2\leqslant
\mathfrak c_{\beta,\ell,\texttt r}m^{-2\beta},
\end{displaymath}
where $p_{\mathcal S_m}^{\perp}$ is the orthogonal projection from $\mathbb L^2(A,dx)$ onto $\mathcal S_m$. If in addition $f_T$ is upper bounded on $A$ by $f_1$, then
\begin{displaymath}
\inf_{\tau\in\mathcal S_m}
\|\tau - b_A\|_{f_T}^{2}
\leqslant f_1\|p_{\mathcal S_m}^{\perp}(b_A) - b_A\|^2
\leqslant\mathfrak c_{\beta,\ell,\texttt r}f_1m^{-2\beta}.
\end{displaymath}
As a consequence, the inequality of Theorem \ref{empirical_risk_bound} can be written
\begin{displaymath}
\mathbb E(\|\widetilde b_m - b_A\|_{N}^{2})
\leqslant\mathfrak c_{\beta,\ell,\texttt r}f_1m^{-2\beta}
+\mathfrak c_{2,H,\sigma}\sigma^2
\frac{\overline{\mathfrak m}_{2,H,M}(T)}{NT^{2 - 2H}}\cdot
\frac{m}{f_0}[1 +\rho(\ell, {\texttt r})m^2] +
\mathfrak c_{\rho,\kappa,\sigma}(1 +\mathfrak b_T)
\frac{\overline{\mathfrak m}_{2,H,M}(T)}{NT}.
\end{displaymath}
We obtain the following result.
%


%
\begin{proposition}\label{empirical_risk_bound_Fourier}
Under Assumption \ref{assumption_b}, if $f_0\leqslant f_T(x)\leqslant f_1$ for every $x \in [\ell, \normalfont{\texttt r}]$, $b_A\in W_{2}^{\beta}([\ell,\normalfont{\texttt r}])$ and $\widetilde b_m$ is computed in the trigonometric basis on $[\ell,\normalfont{\texttt r}]$, then there exists a deterministic constant $\mathfrak c_{\beta,\ell,\normalfont{\texttt r},f_0,f_1} > 0$, not depending on $N$ and $T$, such that with $m = m_{\normalfont{\textrm{opt}}} := [(NT^{2 - 2H})^{1/(2\beta+3)}]$, 
\begin{displaymath}
\mathbb E(\|\widetilde b_{m_{{\rm opt}}} - b_A\|_{N}^{2})
\leqslant
\mathfrak c_{\beta,\ell,\normalfont{\texttt r},f_0,f_1}
\overline{\mathfrak m}_{2,H,M}(T)(NT^{2 - 2H})^{-2\beta /(2\beta + 3)} + 
\mathfrak c_{\rho,\kappa,\sigma}(1 +\mathfrak b_T)
\frac{\overline{\mathfrak m}_{2,H,M}(T)}{NT}.
\end{displaymath}
\end{proposition}
\noindent
We obtain the convergence of the oracle with respect to the empirical risk for a fixed $T$ and $N\rightarrow\infty$, and a rate of convergence which degrades from the rate $N^{-2\beta/(2\beta+1)}$ found in Comte and Genon-Catalot \cite{CGC19} for $H = 1/2$ and $\sigma$ constant, to the rate $N^{ -2\beta /(2\beta+3)}$.
\\
The choice of $m_{{\rm opt}}$ above has the interest to provide a rate, but it  is not possible in practice, as it depends on $\beta$ which is unknown.
\\
\\
Finally, note that the function $b : x\mapsto\mu x$ with $\mu\in\mathbb R^*$ fulfills the conditions of Proposition \ref{empirical_risk_bound_Fourier}. Indeed, since $b' =\mu$, the function $b$ satisfies Assumption \ref{assumption_b} with $m = M =\mu$, and $b_A\in W_{2}^{1}([\ell,\texttt r])$ for every $\texttt r >\ell$. Moreover, since the solution of Equation (\ref{main_equation}) in this case is the fractional Ornstein-Uhlenbeck process, which is a Gaussian process, then for every $\texttt r >\ell$, there exist $f_0,f_1 > 0$ such that $f_0\leqslant f_T\leqslant f_1$. In fact, under Assumption \ref{assumption_b}, thanks to Inequality (\ref{control_density_solution}), $f_T$ is still upper-bounded for nonlinear drift functions.
%


%
\subsubsection{Discussion on the Hermite example}\label{subsubsection_discussion_Hermite}
The second example is the non-compactly supported Hermite basis. Here, $A =\mathbb R$, and the Hermite polynomial and the Hermite function of order $j\geqslant 0$ are given by
\begin{equation}\label{fhermite}
H_j(x) := (-1)^je^{x^2}\frac{d^j}{dx^j}(e^{-x^2})
\textrm{ and }
h_j(x) :=\mathfrak c_jH_j(x)e^{-x^2/2}
\textrm{ $;$ }
\forall x\in\mathbb R,
\end{equation}
where $\mathfrak c_j = (2^jj!\sqrt{\pi})^{-1/2}$.
\\
\\
The sequence $(h_j)_{j\geqslant 0}$ is an orthonormal basis of $\mathbb L^2(\mathbb R,dx)$. By Abramowitz and Stegun \cite{AS64}, and Indritz \cite{INDRITZ61},
\begin{displaymath}
\|h_j\|_{\infty}\leqslant\Phi_0
\textrm{ with }\Phi_0 = 1/\pi^{1/4}.
\end{displaymath} 
So that, for $\varphi_j = h_j$, $L(m)\leqslant\Phi_{0}^{2}m$. In fact, we can prove that $L(m)\leqslant K\sqrt m$ for a constant $K > 0$ (see Comte and Lacour \cite{CL2021}). Moreover, as
\begin{displaymath}
h_j'(x) =\sqrt{\frac{j}{2}}h_{j - 1}(x) -
\frac{\sqrt{j + 1}}2 h_{j + 1}(x),
\end{displaymath}
we find
\begin{displaymath}
\sup_{x\in\mathbb R}
\sum_{j = 0}^{m - 1}
h_j'(x)^2\leqslant 2K m^{3/2}.
\end{displaymath}
Thus, $R(m)\leqslant 2K m^{3/2}$. Here, the bound of Proposition \ref{bound_trace_term} becomes
\begin{displaymath}
\frac{1}{NT}
\textrm{trace}(\mathbf\Psi(m)^{-1}\mathbf\Psi(m,\sigma))
\leqslant
\mathfrak c_{2,H,\sigma}\sigma^2
\frac{\overline{\mathfrak m}_{2,H,M}(T)}{N T^{2 - 2H}}K\sqrt{m}(1+ 2m)\|\mathbf\Psi(m)^{-1}\|_{\textrm{op}},
\end{displaymath}
where $\mathfrak c > 0$ is a universal constant.
\\
\\
This case is more complicated since $f_T$ can no longer be assumed  lower bounded on ${\mathbb R}$, otherwise it would not be integrable.
Therefore, the order of the variance and specifically of $\|\mathbf\Psi(m)^{-1}\|_{\textrm{op}}$ is more difficult to evaluate in general contexts. What is known is that it is growing with $m$ and with order larger than order $\sqrt m$ (see \cite{CGC18}). However, we can still assume that $f_T$ is upper-bounded by a constant $f_1 > 0$, and thus, we can evaluate the bias in a similar way as previously by considering Sobolev-Hermite spaces (see Bongioanni and Torrea \cite{BT06} or Belomestny et al. \cite{BCGC19}) and balls. The Sobolev-Hermite space with regularity $s > 0$ is given by
\begin{equation}\label{Sobolev_Hermite_space}
W_{H}^{s} :=
\left\{\theta\in\mathbb L^2(\mathbb R) :
\sum_{k\geqslant 0}k^sa_k(\theta)^2 <\infty\right\},
\end{equation} 
where $a_k(\theta) :=\langle\theta,h_k\rangle$, $k\in\mathbb N$. The Sobolev-Hermite ball is given by
\begin{displaymath}
W_{H}^{s}(D) :=
\left\{\theta\in\mathbb L^2(\mathbb R) :
\sum_{k\geqslant 0}k^sa_k(\theta)^2\leqslant D\right\}
\textrm{ $;$ }D > 0.
\end{displaymath}
Thus, if $b$ belongs to $W_H^s(D)$, then we have
\begin{displaymath}
\|p_{\mathcal S_m}^{\perp}(b) - b\|^2\leqslant D m^{-s}.
\end{displaymath}
For details, and especially for regularity properties of functions in these spaces, we refer to Section 4.1 of Belomestny et al. \cite{BCGC19}.
%


%
\begin{proposition}\label{empirical_risk_bound_Hermite}
Under Assumption \ref{assumption_b}, if $f_T(x)\leqslant f_1$ for every $x\in\mathbb R$, $b\in W_{H}^{s}(D)$, $\widetilde b_m$ is computed in the Hermite basis on $\mathbb R$, and $\|\mathbf\Psi(m)^{-1}\|_{\normalfont{\textrm{op}}}^2\leqslant m^{\kappa}$ with $\kappa\geqslant 1$, then there exists a deterministic constant $\mathfrak c_{s,D,f_1} > 0$, not depending on $N$ and $T$, such that with $m = m_{\normalfont{\textrm{opt}}} := [(NT^{2 - 2H})^{1/(s + 3/2 +\kappa/2)}]$,
\begin{displaymath}
\mathbb E(\|\widetilde b_{m_{{\rm opt}}} - b\|_{N}^{2}) \leqslant 
\mathfrak c_{s,D,f_1} 
\overline{\mathfrak m}_{2,H,M}(T)(NT^{2 - 2H})^{-2s/(2s + 3 + \kappa)} +
\mathfrak c_{\rho,\kappa,\sigma}(1 +\mathfrak b_T)
\frac{\overline{\mathfrak m}_{2,H,M}(T)}{NT}.
\end{displaymath}
\end{proposition}
\noindent
The rate depends on the unknown $\kappa$, and we mention that if $\|\mathbf\Psi(m)^{-1}\|_{\textrm{op}}$ grows exponentially with $m$, then the rate will become logarithmic, except if the bias also decreases exponentially.
%


%
\subsection{An approximate estimator}\label{subsection_calculable_estimator}
In this subsection, assume that the process $X$ has been observed $N$ times for two close initial conditions $x_0$ and $x_0 +\varepsilon$ with $\varepsilon > 0$, a situation which can occur in pharmacology. We first present a possible application context, and then build an approximate estimator.
%


%
\subsubsection{Application context}\label{subsubsection_application_context}
Let us give details about the application field we have in mind. If $t\mapsto X_{x_0}(t)$ denotes the concentration of a drug along time during its elimination by a patient with initial dose $x_0 > 0$, $t\mapsto X_{x_0 +\varepsilon}(t)$ could be approximated by replicating the exact same protocol on the same patient, but with the initial dose $x_0 +\varepsilon$ after the complete elimination of the previous dose. 
\\
This is an interesting perspective because differential equations driven by the fractional Brownian motion with $H > 1/2$ are well adapted to model the concentration process in pharmacokinetics. Indeed, D'Argenio and Park \cite{DP97} showed that the elimination process has both a deterministic and a random components. A natural way to take into account these two components is to add a stochastic noise in the linear differential equation which classically models the concentration. It has been studied in the It\^o calculus framework by many authors (see e.g. Donnet and Samson \cite{DS13}). However, as mentioned in Delattre and Lavielle \cite{DL11}, the extension of the deterministic concentration model as a diffusion process is not realistic on the biological side because its paths are too rough. 
\\
So, as mentioned in Marie \cite{MARIE14}, Section 5, a way to increase the regularity of the paths of the concentration process is to replace the Brownian motion by a fractional Brownian motion of Hurst index close to $1$ as driving signal.
%


%
\subsubsection{Risk bounds on the approximate estimator}\label{subsubsection_risk_bounds_calculable_estimator}
Throughout this subsection, $b$ fulfills the following reinforced assumption.
%


%
\begin{assumption}\label{assumption_b_reinforced}
The function $b$ belongs to $\normalfont{\textrm{Lip}}_{\mathbf b}^{2}(\mathbb R)$ and fullfils Assumption \ref{assumption_b}.
\end{assumption}
\noindent
Under Assumption \ref{assumption_b_reinforced}, the following proposition allows to approximate the Skorokhod integral of the solution $X$ to Equation (\ref{main_equation}) with respect to $B$ if two paths of $X$ can be observed with different but close initial conditions.
%


%
\begin{proposition}\label{approximation_Skorokhod}
Under Assumption \ref{assumption_b_reinforced}, for every $\varphi\in\normalfont{\textrm{Lip}}_{\mathbf b}^{1}(\mathbb R)$, $\varepsilon > 0$ and $t\in ]0,T]$,
\begin{displaymath}
\left|\int_{0}^{t}\varphi(X_{x_0}(u))\delta X_{x_0}(u) -
S_{\varphi}(x_0,\varepsilon,t)\right|
\leqslant
\alpha_H\sigma^2\frac{\|b''\|_{\infty}\|\varphi'\|_{\infty}}{2}\varepsilon t^{2H - 1}\mathfrak m_{H,M}(t)
\end{displaymath}
where
\begin{eqnarray*}
 S_{\varphi}(x_0,\varepsilon,t) & := &
 \int_{0}^{t}\varphi(X_{x_0}(u))dX_{x_0}(u)\\
 & &
 -\alpha_H\sigma^2
 \int_{0}^{t}
 \int_{0}^{u}\varphi'(X_{x_0}(u))\frac{X_{x_0 +\varepsilon}(u) - X_{x_0}(u)}{X_{x_0 +\varepsilon}(v) - X_{x_0}(v)}|u - v|^{2H - 2}dvdu,
\end{eqnarray*}
$X_x$ is the solution to Equation (\ref{main_equation}) with initial condition $x\in\mathbb R$, and
\begin{displaymath}
\mathfrak m_{H,M}(t) :=
\frac{1}{M^2(2H - 1)}\mathbf 1_{M < 0} +
\frac{t^2}{2H(2H + 1)}\mathbf 1_{M = 0} +
\frac{e^{2Mt}}{M^2(2H - 1)}\mathbf 1_{M > 0}.
\end{displaymath}
\end{proposition}
\noindent
Note that if $M < 0$, then Theorem \ref{approximation_Skorokhod} has been already established in Comte and Marie \cite{CM19} (see Corollary 2.8).
\\
\\
By Proposition \ref{approximation_Skorokhod}, for every $j\in\mathbb N$ and $i\in\{1,\dots,N\}$,
\begin{eqnarray*}
 S_{\varphi_j}^{i}(x_0,\varepsilon,T) & := &
 \int_{0}^{T}\varphi_j(X_{x_0}^{i}(u))dX_{x_0}^{i}(u)\\
 & &
 -\alpha_H\sigma^2
 \int_{0}^{T}
 \int_{0}^{u}\varphi_j'(X_{x_0}^{i}(u))\frac{X_{x_0 +\varepsilon}^{i}(u) - X_{x_0}^{i}(u)}{X_{x_0 +\varepsilon}^{i}(v) - X_{x_0}^{i}(v)}|u - v|^{2H - 2}dvdu
\end{eqnarray*}
provides a good approximation of
\begin{displaymath}
\int_{0}^{T}\varphi_j(X_{x_0}^{i}(s))\delta X_{x_0}^{i}(s).
\end{displaymath}
This legitimates to consider the approximate estimator
\begin{displaymath}
\widehat b_{m,\varepsilon} :=
\sum_{j = 0}^{m - 1}\widehat\theta_{j,\varepsilon}\varphi_j,
\end{displaymath}
where
\begin{displaymath}
\widehat\theta(m,\varepsilon) :=
(\widehat\theta_{0,\varepsilon},\dots,\widehat\theta_{m - 1,\varepsilon})^* =
\widehat{\mathbf\Psi}(m)^{-1}\widehat{\mathbf x}(m,\varepsilon)
\end{displaymath}
and
\begin{displaymath}
\widehat{\mathbf x}(m,\varepsilon) :=
\left(\frac{1}{NT}
\sum_{i = 1}^{N}
S_{\varphi_j}^{i}(x_0,\varepsilon,T)\right)_{j = 0,\dots,m - 1}^{*}.
\end{displaymath}
Contrary to the oracle $\widehat b_m$, $\widehat b_{m,\varepsilon}$ is computable.
%


%
\begin{remark}\label{remark_computation_estimator}
For any $i\in\{1,\dots,N\}$, since $\mathcal I(.,B^i)$ is locally Lipschitz continuous on $\mathbb R$ and $X_{x_0}^{i}$ is $\alpha$-H\"older continuous on $[0,T]$ with $\alpha\in ]1/2,H[$, there exists a square integrable random variable $C :\Omega\rightarrow ]0,\infty[$ such that for any $\eta,\varepsilon\in ]0,1[$,
\begin{eqnarray*}
 |X_{x_0}^{i}(t +\eta) - X_{x_0 +\varepsilon}^{i}(t)|
 & \leqslant &
 |X_{x_0}^{i}(t +\eta) - X_{x_0}^{i}(t)| +
 |X_{x_0}^{i}(t) - X_{x_0 +\varepsilon}^{i}(t)|\\
 & \leqslant &
 C(\eta^{\alpha} +\varepsilon).
\end{eqnarray*}
By taking $\eta =\varepsilon^{1/\alpha}$,
\begin{displaymath}
 |X_{x_0}^{i}(t +\varepsilon^{1/\alpha}) - X_{x_0 +\varepsilon}^{i}(t)|
 \leqslant
 2C\varepsilon.
\end{displaymath}
So, despite the lack of information available on the behavior of the quotient involved in $S_{\varphi_j}^{i}(x_0,\varepsilon,T)$, one could replace it by
\begin{eqnarray*}
 \widetilde S_{\varphi_j}^{i}(x_0,\varepsilon,T) & := &
 \int_{0}^{T}\varphi_j(X_{x_0}^{i}(u))dX_{x_0}^{i}(u)\\
 & &
 -\alpha_H\sigma^2
 \int_{0}^{T}
 \int_{0}^{u}\varphi_j'(X_{x_0}^{i}(u))\frac{X_{x_0}^{i}(u +\varepsilon^{1/\alpha}) - X_{x_0}^{i}(u)}{X_{x_0}^{i}(v +\varepsilon^{1/\alpha}) - X_{x_0}^{i}(v)}|u - v|^{2H - 2}dvdu
\end{eqnarray*}
in the expression of $\widehat\theta(m,\varepsilon)$. This would avoid the requirement of the paths $X_{x_0+\varepsilon}^{i}$ and $X_{x_0}^{i}$ for each individual $i$.
\end{remark}
\noindent
Thus, to be coherent and realistic, we consider the estimator  modified as follows: 
\begin{displaymath}
\widetilde b_{m,\varepsilon} :=
\widehat b_{m,\varepsilon}\mathbf 1_{\widehat\Lambda_{\kappa}(m)}.
\end{displaymath}
Then, we can prove the following result as a consequence of Corollary \ref{risk_bound}.
%


%
\begin{corollary}\label{risk_bound_computable_estimator}
Under Assumptions \ref{assumption_b_reinforced}, \ref{assumption_basis_1} and \ref{assumption_basis_2},
\begin{eqnarray*}
 \mathbb E(\|\widetilde b_{m,\varepsilon} - b_A\|_{f_T}^{2})
 & \leqslant & \mathfrak d_{\kappa, NT}
 \inf_{\tau\in\mathcal S_m}
 \|\tau - b_A\|_{f_T}^{2}
 +\frac{16}{NT}\normalfont{\textrm{trace}}(\mathbf\Psi(m)^{-1}\mathbf\Psi(m,\sigma)) +
 \overline{\mathfrak c}_{\rho,\kappa,\sigma}(1 +\mathfrak b_T)\frac{\overline{\mathfrak m}_{2,H,M}(T)}{NT}\\
 & & 
 +\mathfrak c_{\rho,\kappa,\sigma,H,\mathfrak b''}
 L(m)^{\kappa - 1}\left(\frac{NT}{\log(NT)}\right)^2\frac{\mathfrak m_{H,M}(T)^2}{T^{4 - 4H}}\varepsilon^2
\end{eqnarray*}
where $\mathfrak c_{\rho,\kappa,\sigma,H,\mathfrak b''} > 0$ is a deterministic constant depending only on $\rho$, $\kappa$, $\sigma$, $H$ and $\mathfrak b'' :=\|b''\|_{\infty}^{2}$, and
\begin{displaymath}
\mathfrak d_{\kappa, NT} :=
2\left(1 + 4\frac{3\log(3/2)-1}{(7 +\kappa)\log(NT)}\right).
\end{displaymath}
\end{corollary}
\noindent
In order to keep the rate of convergence obtained for $\widetilde b_m$, one can assume that $\varepsilon$ depends on $N$ and $T$, and take
\begin{displaymath}
\varepsilon_{N,T} =
\frac{1}{(NT)^{1/2}}
\left(\frac{NT}{\log(NT)}\right)^{-1}.
\end{displaymath}
%


%
\section{An alternative estimator}\label{section_alternative_estimator}
The cornerstone of the proof of Proposition \ref{approximation_Skorokhod} is that, for every $t\in [0,T]$,
\begin{displaymath}
Y_x(t) =
\int_{0}^{t}b'(X_x(s))ds
\textrm{ with }
Y_x(t) :=
\log(\partial_xX_x(t)).
\end{displaymath}
By assuming in the first place that $X_x$ and $\partial_xX_x$ can be observed $N$ times on $[0,T]$, the previous equality suggests the noiseless model
\begin{equation}\label{alternative_model}
Y_{x}^{i}(T) =
\int_{0}^{T}b'(X_{x}^{i}(s))ds
\textrm{ $;$ }i\in\{1,\dots,N\},
\end{equation}
where $X_{x}^{1},\dots,X_{x}^{N}$ (resp. $Y_{x}^{1},\dots,Y_{x}^{N}$) are $N$ independent copies of $X_x$ (resp. $Y_x$).
\\
\\
Consider the objective function $\gamma_{N}^{\dag}$ defined by
\begin{displaymath}
\gamma_{N}^{\dag}(\tau) :=
\frac{1}{NT}\sum_{i = 1}^{N}\left(
\int_{0}^{T}\tau(X_{x}^{i}(s))^2ds - 2\int_{0}^{T}\tau(X_{x}^{i}(s))dY_{x}^{i}(s)\right)
\end{displaymath}
for every function $\tau :\mathbb R\rightarrow\mathbb R$.
\\
\\
Note that for any bounded function $\tau$ from $\mathbb R$ into itself,
\begin{eqnarray*}
 \mathbb E(\gamma_{N}^{\dag}(\tau)) & = &
 \frac{1}{T}\int_{0}^{T}\mathbb E(\tau(X_x(s))^2)ds -
 \frac{2}{T}\int_{0}^{T}\mathbb E(\tau(X_x(s))b'(X_x(s)))ds\\
 & = &
 \frac{1}{T}\int_{0}^{T}\mathbb E((\tau(X_x(s)) - b'(X_x(s)))^2)ds -\frac{1}{T}\int_{0}^{T}\mathbb E(b'(X_x(s))^2)ds.
\end{eqnarray*}
Then, the definition of $f_T$ gives
\begin{equation}\label{EspGamma_2}
\mathbb E(\gamma_{N}^{\dag}(\tau)) =
\int_{-\infty}^{\infty}(\tau(x) - b'(x))^2f_T(x)dx -
\int_{-\infty}^{\infty}b'(x)^2f_T(x)dx.
\end{equation}
Equality (\ref{EspGamma_2}) shows that $\mathbb E(\gamma_{N}^{\dag}(\tau))$ is the smallest for $\tau$ the nearest of $b'$. Therefore, to minimize its empirical version $\gamma_{N}^{\dag}(\tau)$ should provide a functional estimator near of $b'$.
\\
\\
Let
\begin{displaymath}
\widehat b_{m}^{\dag,\prime} =\arg\min_{\tau\in\mathcal S_m}\gamma_{N}^{\dag}(\tau)
\end{displaymath}
be the projection estimator of $b_A' := b_{|A}'$ on $\mathcal S_m$. Then,
\begin{displaymath}
\widehat b_{m}^{\dag,\prime} =
\sum_{j = 0}^{m - 1}\widehat\theta_j'\varphi_j
\end{displaymath}
where
\begin{displaymath}
\widehat{\theta}'(m) :=
(\widehat\theta_0',\dots,\widehat\theta_{m - 1}')^* =
\widehat{\mathbf\Psi}(m)^{-1}\widehat{\mathbf y}(m)
\end{displaymath}
with
\begin{displaymath}
\widehat{\mathbf y}(m) :=
\left(\frac{1}{NT}
\sum_{i = 1}^{N}
\int_{0}^{T}\varphi_j(X_{x}^{i}(s))dY_{x}^{i}(s)\right)_{j = 0,\dots,m - 1}^{*} =
(\langle\varphi_j,b'\rangle_N)_{j = 0,\dots,m - 1}^{*}.
\end{displaymath}
As for $\widehat b_m$, in order to ensure the existence and the stability of the estimator, $\widehat b_{m}^{\dag,\prime}$ is replaced by
\begin{displaymath}
\widetilde b_{m}^{\dag,\prime} :=
\widehat b_{m}^{\dag,\prime}\mathbf 1_{\widehat\Lambda_0(m)},
\end{displaymath}
where $\widehat\Lambda_0(m)$ is defined by (\ref{Lambdakappa}) with $\kappa = 0$, and $L(m)$ fulfills the following assumption. In the continuous-time processes framework, this estimator is similar to the least square estimator for the noiseless regression model of Cohen et al. \cite{CDL13}.
%


%
\begin{assumption}\label{assumption_basis_2_dag}
The quantity $L(m)$ satisfies
\begin{displaymath}
L(m)
(\|\mathbf\Psi(m)^{-1}\|_{\normalfont\textrm{op}}\vee 1)
\leqslant
\frac{\mathfrak c_{0,T}}{2}\cdot\frac{NT}{\log(NT)}
\quad
\textrm{with}
\quad
\mathfrak c_{0,T} =\frac{3\log(3/2)-1}{7T}.
\end{displaymath}
\end{assumption}
\noindent
The following result provides controls of the empirical risk and of the $f_T$-weighted integrated  risk of $\widetilde b_{m}^{\dag,\prime}$ respectively.
%


%
\begin{proposition}\label{empirical_risk_bound_estimator_dag}
Under Assumptions \ref{assumption_b}, \ref{assumption_basis_1} and \ref{assumption_basis_2_dag},
\begin{equation}\label{bounddag1}
\mathbb E(\|\widetilde b_{m}^{\dag,\prime} - b_A'\|_{N}^{2})
\leqslant
\inf_{\tau\in\mathcal S_m}
\|\tau - b_A'\|_{f_T}^{2} +
\frac{\mathfrak c_{\ref{empirical_risk_bound_estimator_dag}}}{NT}
\end{equation}
where $\mathfrak c_{\ref{empirical_risk_bound_estimator_dag}} :=\mathfrak c_{\Omega}^{1/2}\|b'\|_{\infty}^{2}$. In addition,
\begin{equation}\label{bounddag2}
\mathbb E(\|\widetilde b_{m}^{\dag,\prime} - b_A'\|_{f_T}^{2})
\leqslant
\left(1 +  \frac{2T\mathfrak c_{0,T}}{\log(NT)}\right)
\inf_{\tau\in\mathcal S_m}
\|\tau - b_A'\|_{f_T}^{2}
+\overline{\mathfrak c}_{\mathfrak b'}\frac{1 +\mathfrak c_{0,T}}{NT}
\end{equation}
where $\overline{\mathfrak c}_{\mathfrak b'} > 0$ is a deterministic constant depending only on $\mathfrak b' :=\|b'\|_{\infty}$.
\end{proposition}
\noindent
A consequence of Proposition \ref{empirical_risk_bound_estimator_dag} is that the estimator can reach a parametric rate and a risk of order $1/(NT)$. No model selection step is required, one has only to choose the largest $m$ of the collection. In practice, this leads to a good but non parsimonious estimator. A theoretical stopping rule leading to keep only an adequate number of projection coefficients would be interesting but is not available yet.
\\
\\
Now, assume that the process $X$ has been observed $N$ times for two close initial conditions $x_0$ and $x_0 +\varepsilon$ with $\varepsilon > 0$. Let us consider the estimator
\begin{displaymath}
\widehat b_{m,\varepsilon}^{\dag,\prime} :=
\sum_{j = 0}^{m - 1}\widehat\theta_{j,\varepsilon}'\varphi_j
\end{displaymath}
with
\begin{displaymath}
\widehat\theta'(m,\varepsilon) :=
(\widehat\theta_{0,\varepsilon}',\dots,\widehat\theta_{m - 1,\varepsilon}') =
\widehat{\mathbf\Psi}(m)^{-1}\widehat{\mathbf y}(m,\varepsilon)
\end{displaymath}
and
\begin{displaymath}
\widehat{\mathbf y}(m,\varepsilon) :=
\left(
\frac{1}{NT}
\sum_{i = 1}^{N}
\int_{0}^{T}\varphi_j(X_{x_0}^{i}(s))dY_{x_0,\varepsilon}^{i}(s)
\right)_{j = 0,\dots,m - 1}^{*},
\end{displaymath}
where $Y_{x_0,\varepsilon}^{1},\dots,Y_{x_0,\varepsilon}^{N}$ are $N$ independent copies of the process $Y_{x_0,\varepsilon}$ defined by
\begin{displaymath}
Y_{x_0,\varepsilon}(t) :=
\log\left(
\frac{X_{x_0 +\varepsilon}(t) - X_{x_0}(t)}{\varepsilon}\right)
\textrm{ $;$ }
\forall t\in [0,T].
\end{displaymath}
%


%
\begin{corollary}\label{risk_bound_computable_estimator_dag}
Under Assumptions \ref{assumption_b_reinforced}, \ref{assumption_basis_1} and \ref{assumption_basis_2_dag},
\begin{eqnarray*}
 \mathbb E(\|\widetilde b_{m,\varepsilon}^{\dag,\prime} - b_A'\|_{f_T}^{2})
 & \leqslant & \mathfrak d_{0, NT}
 \inf_{\tau\in\mathcal S_m}
 \|\tau - b_A'\|_{f_T}^{2} +
 \overline{\mathfrak c}_{\mathfrak b'}\frac{1 +\mathfrak c_{0,T}}{NT}\\
 & & 
 \hspace{2cm}
 +\frac{\|b''\|_{\infty}^{2}}{2}\mathfrak c_{0,T}^{2}e^{2MT}\frac{1}{L(m)}\left(\frac{NT}{\log(NT)}\right)^2\varepsilon^2.
\end{eqnarray*}
\end{corollary}
\noindent
Finally, we have obtained a good estimator of $b'_A$, but there remains difficulties to deduce an estimator of $b_A$. Assume that $A = [\texttt r,\ell]$ with $\texttt r,\ell\in\mathbb R$ and $\texttt r <\ell$, and that there exists $f_0 > 0$ such that $f_T(x)\geqslant f_0$ for every $x\in A$. For every $x\in A$, consider the estimator
\begin{displaymath}
\widehat b_{m,\varepsilon}^{\dag}(x) := \widehat b_A(\texttt r) +
\int_{\texttt r}^{x}\widehat b_{m,\varepsilon}^{\dag,\prime}(y)dy
\end{displaymath}
of
\begin{displaymath}
b_A(x) := b_A(\texttt r) +\int_{\texttt r}^{x}b_A'(y)dy
\textrm{ $;$ }
x\in A.
\end{displaymath}
Then,
\begin{eqnarray*}
 \|\widehat b_{m,\varepsilon}^{\dag} - b_A\|_{f_T}^{2} & = &
 \int_{\texttt r}^{\ell}\left[(\widehat b_A(\texttt r) - b_A(\texttt r)) + \int_{\texttt r}^{x}(\widehat b_{m,\varepsilon}^{\dag,\prime}(y) - b_A'(y))dy\right]^2 f_T(x)dx\\
 & \leqslant & 2(\widehat b_A(\texttt r)-b_A(\texttt r))^2 + 2(\ell -\texttt r)\left(\int_{\texttt r}^{\ell}f_T(x)dx\right)\left(\int_{\texttt r}^{\ell}|\widehat b_{m,\varepsilon}^{\dag,\prime}(y) - b_A'(y)|^2dy\right) \\
 & \leqslant & 2(\widehat b_A(\texttt r)-b_A(\texttt r))^2 + 2\frac{\ell -\texttt r}{f_0}\|\widehat b_{m,\varepsilon}^{\dag,\prime} - b_A'\|_{f_T}^{2}.
\end{eqnarray*}
So, thanks to Corollary \ref{risk_bound_computable_estimator_dag},
\begin{eqnarray*}
 \mathbb E(\|\widetilde b_{m,\varepsilon}^{\dag} - b_A\|_{f_T}^{2})
 & \leqslant & 2(\widehat b_A(\texttt r)-b_A(\texttt r))^2 + 2\frac{\ell -\texttt r}{f_0}\left[
 \mathfrak d_{0, NT}
 \inf_{\tau\in\mathcal S_m}
 \|\tau - b_A'\|_{f_T}^{2} +
 \overline{\mathfrak c}_{\mathfrak b'}\frac{1 +\mathfrak c_{0,T}}{NT}\right.\\
 & & 
 \hspace{2cm}
 \left.
 +\frac{\|b''\|_{\infty}^{2}}{2}\mathfrak c_{0,T}^{2}e^{2MT}\frac{1}{L(m)}\left(\frac{NT}{\log(NT)}\right)^2\varepsilon^2\right].
\end{eqnarray*}
Clearly, $\widetilde b_{m,\varepsilon}^{\dag}$ keeps a good rate of convergence only if the function $b_A$ is known at one point, or if $b_A$ can be estimated at one point with a good rate. No pointwise procedure which would grant to keep the parametric feature is available. Note also that, even if $b_A(\texttt r)$ is known, another drawback of $\widetilde b_{m,\varepsilon}^{\dag}$ is that, since it is obtained by integrating an estimator of $b_A'$, there are restrictive conditions: $A$ must be a compact interval with $f_T\geqslant f_0 > 0$ on $A$.
%


%
\section{Concluding remarks}\label{concluding}
In this paper, we propose two estimation strategies for nonparametric reconstruction of the drift function $b$ from $N$ i.i.d. continuous observations drawn in the fractional SDE given by (\ref{main_equation}). On the one hand, we bound the empirical and $f_T$-weighted $\mathbb L^2$-risks of the auxiliary oracle $\widehat b_m$, and then of the approximate estimator $\widehat b_{m,\varepsilon}$. We compare the result with non fractional case. In the case of a specific trigonometric basis, we can, under additional assumptions, prove the consistency of the estimator and evaluate its rate of convergence. However, the choice of $m_{\textrm{opt}}$ which is proposed to obtain this result is not possible in practice, as it depends on unknown parameters. Therefore, a model selection step in the spirit of Comte and Genon-Catalot \cite{CGC19} would have to be settled. On the other hand, we bound the empirical and $f_T$-weighted $\mathbb L^2$-risks of the alternative calculable estimator $\widehat b_{m,\varepsilon}^{\dag,\prime}$, get a parametric rate of convergence, and show that we can not deduce  a $f_T$-weighted $\mathbb L^2$-risk bound on the primitive estimator $\widehat b_{m,\varepsilon}^{\dag}$, except if the function is known at one point. Even so, there will be the restrictive conditions requiring that $A$ is a compact interval and that $f_T$ is lower-bounded by a strictly positive constant on $A$. For this reason, $\widehat b_{m,\varepsilon}$ is adequate to estimate $b$ and the other strategy to estimate its derivative, which has interest also in many applications.
\\
For both estimators $\widehat b_{m,\varepsilon}$ and $\widehat b_{m,\varepsilon}^{\dag}$, two paths with slightly different initial conditions have to be available for each individual $i$, but this context may be realistic in pharmacokinetics experiments. We also propose another idea of transformation of $\widehat b_m$ into a computable version, which seems intuitive and does not require two paths per individual, but would require a deeper study (see Remark \ref{remark_computation_estimator}).
\\
Lastly, the tedious question of discretization may be investigated in the future, to take into account the fact that, for each individual $i$, the observation is $(X_i(k\Delta))_{1\leqslant k\leqslant n}$ where $n\Delta = T$.
%
 

%
\section{Proofs}\label{section_proofs}
%


%
\subsection{Proof of Theorem \ref{control_divergence_integral}}
On the one hand, for any $s,t\in [0,T]$,
\begin{displaymath}
\mathbf D_sX(t) =
\sigma\mathbf 1_{[0,t]}(s)\exp\left(\int_{s}^{t}b'(X(u))du\right)
\end{displaymath}
and, by Assumption \ref{assumption_b},
\begin{displaymath}
|\mathbf D_sX(t)|
\leqslant
\sigma
\mathbf 1_{[0,t]}(s)
e^{M(t - s)}.
\end{displaymath}
Then, by the chain rule for Malliavin's derivative and Jensen's inequality,
\begin{eqnarray}
 \mathbb E\left(\int_{0}^{T}
 \int_{0}^{T}|\mathbf D_s[\varphi(X(t))]|^{1/H}dsdt\right)^{pH}
 & = &
 \mathbb E\left(\int_{0}^{T}
 \int_{0}^{T}
 |\varphi'(X(t))\mathbf D_sX(t)|^{1/H}dsdt\right)^{pH}
 \nonumber\\
 & \leqslant &
 \sigma^p\left(\int_{0}^{T}
 \mathbb E(|\varphi'(X(t))|^{1/H})
 \int_{0}^{t}
 e^{M/H(t - s)}dsdt\right)^{pH}
 \nonumber\\
 \label{control_divergence_integral_1}
 & \leqslant &
 \sigma^p
 \mathfrak m_{p,H,M}(T)\left(\int_{0}^{T}
 \mathbb E(|\varphi'(X(t))|^p)^{{1/(pH)}}dt\right)^{pH}
\end{eqnarray}
where
\begin{displaymath}
\mathfrak m_{p,H,M}(T) =
\left(-\frac{H}{M}\right)^{pH}\mathbf 1_{M < 0} +
T^{pH}\mathbf 1_{M = 0} +
\left(\frac{H}{M}\right)^{pH}e^{pMT}\mathbf 1_{M > 0}.
\end{displaymath}
On the other hand, by Hu et al. \cite{HNZ18}, Lemma 3.1, there exists a deterministic constant $\mathfrak c_{p,H} > 0$, depending only on $p$ and $H$, such that for any $\varphi\in\textrm{Lip}_{\mathbf b}^{1}(\mathbb R)$,
\begin{eqnarray}
 \mathbb E\left(
 \left|\int_{0}^{T}\varphi(X(t))\delta B(t)\right|^p\right)
 & \leqslant &
 \mathfrak c_{p,H}\left[\left(\int_{0}^{T}\mathbb E(|\varphi(X(t))|^{1/H})dt\right)^{pH}\right.
 \nonumber\\
 \label{control_divergence_integral_2}
 & &
 +\left.
 \mathbb E\left(\int_{0}^{T}
 \int_{0}^{T}|\mathbf D_s[\varphi(X(t))]|^{1/H}dsdt\right)^{pH}\right].
\end{eqnarray}
Inequalities (\ref{control_divergence_integral_1}) and (\ref{control_divergence_integral_2}) together allow to conclude.
%


%
\subsection{Proof of Proposition \ref{density_solution_main_equation}}
Consider $t\in ]0,T]$. On the one hand, for any $s\in\mathbb R_+$,
\begin{displaymath}
\mathbf D_sX(t) =
\sigma\mathbf 1_{[0,t]}(s)
\exp\left(\int_{s}^{t}b'(X(u))du\right).
\end{displaymath}
So, since $b'$ is a $[m,M]$-valued function,
\begin{displaymath}
\sigma\mathbf 1_{[0,t]}(s)e^{m(t - s)}
\leqslant\mathbf D_sX(t)
\leqslant
\sigma\mathbf 1_{[0,t]}(s)e^{M(t - s)}.
\end{displaymath}
On the other hand, by Nourdin and Viens \cite{NV09}, Proposition 3.7,
\begin{displaymath}
-\mathbf D_s\mathbf L^{-1}X(t) =
\int_{0}^{\infty}
e^{-u}\mathbf T_u(\mathbf D_sX(t))du
\end{displaymath}
where $(\mathbf T_u)_{u\in\mathbb R_+}$ is the Ornstein-Uhlenbeck semigroup (see Nualart \cite{NUALART06}, Section 1.4). Moreover, for any $u\in\mathbb R_+$, $(\mathbf T_u)_{|\mathbb R} =\textrm{Id}_{\mathbb R}$ by Mehler's formula (see Nualart \cite{NUALART06}, Equation (1.67)), and for every $U_1,U_2\in\mathbb L^2(\Omega,\mathbb P)$,
\begin{displaymath}
U_1\geqslant U_2\Longrightarrow
\mathbf T_u(U_1)\geqslant\mathbf T_u(U_2)
\end{displaymath}
by Nualart \cite{NUALART06}, Property (i) page 55. Then,
\begin{displaymath}
\sigma\mathbf 1_{[0,t]}(s)e^{m(t - s)}
\leqslant
-\mathbf D_s\mathbf L^{-1}X(t)
\leqslant
\sigma\mathbf 1_{[0,t]}(s)e^{M(t - s)}.
\end{displaymath}
Therefore,
\begin{displaymath}
\sigma(m,t)^2
\leqslant
g_{t}^{*}(x_0,X^*(t))
\leqslant
\sigma(M,t)^2
\end{displaymath}
where
\begin{displaymath}
\sigma(\mu,t)^2 =
\alpha_H\sigma^2\int_{0}^{t}\int_{0}^{t}|v - u|^{2H - 2}
e^{\mu(2t - v - u)}dudv > 0
\textrm{ $;$ }
\forall\mu\in\mathbb R.
\end{displaymath}
Nourdin and Viens \cite{NV09}, Theorem 3.1 and Corollary 3.5 allow to conclude.
%


%
\subsection{Proof of Theorem \ref{empirical_risk_bound}}
The proof relies on two lemmas which are stated first.
%


%
\begin{lemma}\label{bound_e}
There exists a deterministic constant $\mathfrak c_{\mathbf e} > 0$, not depending on $m$, $N$ and $T$, such that
\begin{displaymath}
\mathbb E(|\mathbf e(m)^*\mathbf e(m)|^2)
\leqslant
\mathfrak c_{\mathbf e}\sigma^4
\frac{\overline{\mathfrak m}_{4,H,M}(T)}{N^2T^{4 - 4H}}m(L(m)^2 + R(m)^2).
\end{displaymath}
\end{lemma}
%


%
\begin{lemma}\label{bound_Lambda}
Consider
\begin{displaymath}
\Omega(m) :=
\left\{\sup_{\tau\in\mathcal S_m}
\left|\frac{\|\tau\|_{N}^{2}}{\|\tau\|_{f_T}^{2}} - 1\right|\leqslant\frac{1}{2}\right\}.
\end{displaymath}
Under Assumptions \ref{assumption_b}, \ref{assumption_basis_1} and \ref{assumption_basis_2}, there exists a deterministic constant $\mathfrak c_{\Omega} > 0$, not depending on $m$, $N$ and $T$, such that
\begin{displaymath}
\mathbb P(\widehat\Lambda_{\kappa}(m)^c)\leqslant\frac{\mathfrak c_{\Omega}}{(NT)^{6 +\kappa}}
\textrm{ and }
\mathbb P(\Omega(m)^c)\leqslant\frac{\mathfrak c_{\Omega}}{(NT)^{6 +\kappa}}.
\end{displaymath}
\end{lemma}
%


%
\subsubsection{Steps of the proof}
First of all,
\begin{equation}\label{empirical_risk_bound_1}
\|\widetilde b_m - b_A\|_{N}^{2} =
\|\widehat b_m - b_A\|_{N}^{2}\mathbf 1_{\widehat\Lambda_{\kappa}(m)} +
\|b_A\|_{N}^{2}\mathbf 1_{\widehat\Lambda_{\kappa}(m)^c} =
U_1 + U_2 + U_3
\end{equation}
where
\begin{eqnarray*}
 U_1 & := &
 \|b_A\|_{N}^{2}\mathbf 1_{\widehat\Lambda_{\kappa}(m)^c},\\
 U_2 & := &
 \|\widehat b_m - b_A\|_{N}^{2}\mathbf 1_{\widehat\Lambda_{\kappa}(m)\cap\Omega(m)}
 \textrm{ and}\\
 U_3 & := &
 \|\widehat b_m - b_A\|_{N}^{2}\mathbf 1_{\widehat\Lambda_{\kappa}(m)\cap\Omega(m)^c}.
\end{eqnarray*}
Let us find suitable bounds on $\mathbb E(U_1)$, $\mathbb E(U_2)$ and $\mathbb E(U_3)$.
\begin{itemize}
 \item\textbf{Bound for $\mathbb E(U_1)$.} By Lemma \ref{bound_Lambda},
 \begin{displaymath}
 \mathbb E(U_1)
 \leqslant
 \mathbb E(\|b_A\|_{N}^{4})^{1/2}
 \mathbb P(\widehat\Lambda_{\kappa}(m)^c)^{1/2}
 \leqslant
 \frac{\mathfrak c_1}{NT}
 \end{displaymath}
 where
 \begin{displaymath}
 \mathfrak c_1 :=
 \mathfrak c_{\Omega}^{1/2}
 \left(\int_{-\infty}^{\infty}b_A(x)^4f_T(x)dx\right)^{1/2} =
 \mathfrak c_{\Omega}^{1/2}\|b_{A}^{2}\|_{f_T}.
 \end{displaymath}
 \item\textbf{Bound for $\mathbb E(U_2)$.} By denoting $p_{N,\mathcal S_m}^{\perp}$ the orthogonal projection from $\mathbb L^2(A,f_T(x)dx)$ onto $\mathcal S_m$ with respect to the empirical scalar product $\langle .,.\rangle_N$,
 \begin{equation}\label{empirical_risk_bound_2}
 \|\widehat b_m - b_A\|_{N}^{2} =
 \|\widehat b_m - p_{N,\mathcal S_m}^{\perp}(b_A)\|_{N}^{2}
 +\inf_{\tau\in\mathcal S_m}\|b_A -\tau\|_{N}^{2}.
 \end{equation}
 As in the proof of Comte and Genon-Catalot \cite{CGC19}, Proposition 2.1, on $\Omega(m)$,
 \begin{displaymath}
 \|\widehat b_m - p_{N,\mathcal S_m}^{\perp}(b_A)\|_{N}^{2}
 =\mathbf e(m)^*\widehat{\mathbf\Psi}(m)^{-1}\mathbf e(m)
 \leqslant
 2\mathbf e(m)^*\mathbf\Psi(m)^{-1}\mathbf e(m).
 \end{displaymath}
 So,
\begin{eqnarray*}
  \mathbb E(
  \|\widehat b_m - p_{N,\mathcal S_m}^{\perp}(b_A)\|_{N}^{2}
  \mathbf 1_{\Omega(m)\cap\widehat\Lambda_{\kappa}(m)})
  & \leqslant &
  2\mathbb E\left(\sum_{j,k = 0}^{m - 1}
  \mathbf e(m)_j\mathbf e(m)_k
  \mathbf\Psi(m)_{j,k}^{-1}\right)\\
  & = &
  \frac{2\sigma^2}{NT^2}
  \sum_{j,k = 0}^{m - 1}\mathbf\Psi(m)_{j,k}^{-1}\\
  & &
  \times\mathbb E\left(\left(\int_{0}^{T}\varphi_j(X(s))\delta B(s)\right)
  \left(\int_{0}^{T}\varphi_k(X(s))\delta B(s)\right)\right)\\
  & = &
  \frac{2}{NT}
  \textrm{trace}(\mathbf\Psi(m)^{-1}\mathbf\Psi(m,\sigma)).
 \end{eqnarray*}
 Therefore, by Inequality (\ref{empirical_risk_bound_2}),
 \begin{displaymath}
 \mathbb E(U_2)\leqslant
 \inf_{\tau\in\mathcal S_m}\|b_A -\tau\|_{f_T}^{2} +
 \frac{2}{NT}
 \textrm{trace}(\mathbf\Psi(m)^{-1}\mathbf\Psi(m,\sigma)).
 \end{displaymath}
 \item\textbf{Bound for $\mathbb E(U_3)$.} On the one hand,
 \begin{displaymath}
 \|\widehat b_m - p_{N,\mathcal S_m}^{\perp}(b_A)\|_{N}^{2} =
 \mathbf e(m)^*\widehat{\mathbf\Psi}(m)^{-1}\mathbf e(m).
 \end{displaymath}
 On the other hand, on $\widehat\Lambda_{\kappa}(m)$, for $N,T$ large enough,
 \begin{eqnarray*}
  (L(m) + R(m))\|\widehat{\mathbf\Psi}(m)^{-1}\|_{\textrm{op}}
  & \leqslant &
  L(m)(\|\widehat{\mathbf\Psi}(m)^{-1}\|_{\textrm{op}}\vee 1) +\rho L(m)^{\kappa}(\|\widehat{\mathbf\Psi}(m)^{-1}\|_{\textrm{op}}\vee 1)^{\kappa}\\
  & \leqslant &
  (1 +\rho)\mathfrak c_{\kappa,T}\left(\frac{NT}{\log(NT)}\right)^{\kappa}.
 \end{eqnarray*}
 Then, by Lemmas \ref{bound_e} and \ref{bound_Lambda}, there exists a deterministic constant $\mathfrak c_2 > 0$, not depending on $m$ (satisfying $m\leqslant NT$), $N$ and $T$, such that
\begin{eqnarray*}
  \mathbb E(U_3) & \leqslant &
  \mathbb E((
  \|\widehat b_m - p_{N,\mathcal S_m}^{\perp}(b_A)\|_{N}^{2}
  +\|b_A\|_{N}^{2})\mathbf 1_{\widehat\Lambda_{\kappa}(m)\cap\Omega(m)^c})\\
  & \leqslant &
  \left[
  \frac{(1 +\rho)\mathfrak c_{\kappa,T}}{L(m) + R(m)}\left(\frac{NT}{\log(NT)}\right)^{\kappa}
  \mathbb E(|\mathbf e(m)^*\mathbf e(m)|^2)^{1/2} +
  \mathbb E(\|b_A\|_{N}^{4})^{1/2}\right]
  \mathbb P(\Omega(m)^c)^{1/2}\\
  & \leqslant &
  \mathfrak c_2\frac{\overline{\mathfrak m}_{2,H,M}(T)}{NT}.
 \end{eqnarray*}
\end{itemize}
These bounds together with Inequality (\ref{empirical_risk_bound_1}) allow to conclude.
%


%
\subsubsection{Proof of Lemma \ref{bound_e}}
On the one hand, for any bounded and measurable function $\tau :\mathbb R\rightarrow\mathbb R$,
\begin{equation}\label{bound_e_1}
\int_{0}^{T}
\mathbb E(\tau(X(s))ds =
\int_{-\infty}^{\infty}\int_{0}^{T}
\tau(x)p_s(x_0,x)dsdx =
T\int_{-\infty}^{\infty}\tau(x)f_T(x)dx
\end{equation}
and, by Jensen's inequality, for every $p > 1/H$,
\begin{eqnarray}
 \int_{0}^{T}
 \mathbb E(\tau(X(s)))^{1/(pH)}ds & = &
 T\int_{0}^{T}\left(\int_{-\infty}^{\infty}\tau(x)p_s(x_0,x)dx\right)^{1/(pH)}\frac{ds}{T}
 \nonumber\\
 & \leqslant &
 T\left(\int_{-\infty}^{\infty}\int_{0}^{T}
 \tau(x)p_s(x_0,x)\frac{ds}{T}dx\right)^{1/(pH)}
 \nonumber\\
 \label{bound_e_2}
 & = &
 T\left(\int_{-\infty}^{\infty}\tau(x)f_T(x)dx\right)^{1/(pH)}.
\end{eqnarray}
On the other hand, by Jensen's inequality, since $(B_1,X_1),\dots,(B_N,X_N)$ are i.i.d and by Equality (\ref{divergence_operator_1}),
\begin{eqnarray*}
 \mathbb E(|\mathbf e(m)^*\mathbf e(m)|^2) & = &
 \mathbb E\left(\left|\sum_{j = 0}^{m - 1}\mathbf e(m)_{j}^{2}\right|^2\right)\\
 & \leqslant &
 \frac{\sigma^4m}{(NT)^4}\sum_{j = 0}^{m - 1}\mathbb E\left[
 \left(\sum_{i = 1}^{N}\int_{0}^{T}\varphi_j(X^i(s))\delta B^i(s)\right)^4\right]\\
 & \leqslant &
 \frac{\sigma^4m}{N^3T^4}\sum_{j = 0}^{m - 1}\mathbb E(Y_{j}^{4})
 +\frac{\sigma^4m}{N^2T^4}\sum_{j = 0}^{m - 1}\mathbb E(Y_{j}^{2})^2
\end{eqnarray*}
where
\begin{displaymath}
Y_j :=
\int_{0}^{T}\varphi_j(X(s))\delta B(s)
\textrm{ $;$ }
\forall j\in\{0,\dots,m - 1\}.
\end{displaymath}
 Moreover, by Theorem \ref{control_divergence_integral}, Equality (\ref{bound_e_1}) and Inequality (\ref{bound_e_2}), for $j = 0,\dots,m - 1$,
\begin{eqnarray*}
 \mathbb E(Y_{j}^{4})
 & \leqslant &
 \mathfrak c_{4,H,\sigma}\overline{\mathfrak m}_{4,H,M}(T)
 \left[\left(\int_{0}^{T}\mathbb E(|\varphi_j(X(s))|^{1/H})ds\right)^{4H}
 +\left(\int_{0}^{T}\mathbb E(|\varphi_j'(X(s))|^4)^{1/(4H)}ds\right)^{4H}\right]\\
 & \leqslant &
 \mathfrak c_{4,H,\sigma}\overline{\mathfrak m}_{4,H,M}(T)T^{4H}\left[
 \left(\int_{-\infty}^{\infty}|\varphi_j(x)|^{1/H}f_T(x)dx\right)^{4H}
 +\int_{-\infty}^{\infty}|\varphi_j'(x)|^4f_T(x)dx\right]\\
 & \leqslant &
 \mathfrak c_{4,H,\sigma}\overline{\mathfrak m}_{4,H,M}(T)
 T^{4H}\sup_{x\in A}\{|\varphi_j(x)|^4 + |\varphi_j'(x)|^4\}.
\end{eqnarray*}
Therefore, there exists a deterministic constant $\mathfrak c_{\mathbf e} > 0$, not depending on $m$, $N$ and $T$, such that
\begin{displaymath}
\mathbb E(|\mathbf e(m)^*\mathbf e(m)|^2)
\leqslant
\mathfrak c_{\mathbf e}\sigma^4
\frac{\overline{\mathfrak m}_{4,H,M}(T)}{N^2T^{4 - 4H}}m(L(m)^2 + R(m)^2).
\end{displaymath}
%


%
\subsubsection{Proof of Lemma \ref{bound_Lambda}}
The proof of Lemma \ref{bound_Lambda} is close to the proof of Comte and Genon-Catalot \cite{CGC19}, Lemma 3.
\\
\\
On the one hand, by the beginning of Comte and Genon-Catalot \cite{CGC19}, Subsection 6.1 which remains true for $H > 1/2$ without additional arguments,
\begin{displaymath}
\Omega(m) =
\left\{\|\mathbf\Psi(m)^{-1/2}\widehat{\mathbf\Psi}(m)
\mathbf\Psi(m)^{-1/2} -\mathbf I_m\|_{\textrm{op}}\leqslant\frac{1}{2}\right\}.
\end{displaymath}
Then, as in the proof of Comte and Genon-Catalot \cite{CGC18}, Proposition 3, by Cohen et al. \cite{CDL13}, Theorem 1,
\begin{displaymath}
\mathbb P(\Omega(m)^c)
\leqslant
2m\exp\left(
-\mathfrak c_1
\frac{NT}{L(m)(\|\mathbf\Psi(m)^{-1}\|_{\textrm{op}}\vee 1)}
\right)
\end{displaymath}
with $\mathfrak c_1 := 1/2(3\log(3/2)-1)$. Under Assumption \ref{assumption_basis_2},
\begin{displaymath}
L(m)
(\|\mathbf\Psi(m)^{-1}\|_{\normalfont\textrm{op}}\vee 1)
\leqslant
\frac{\mathfrak c_{\kappa,T}}{2}\cdot\frac{NT}{\log(NT)}.
\end{displaymath}
Then, there exists a deterministic constant $\mathfrak c_{\Omega} > 0$, not depending on $m$, $N$ and $T$, such that
\begin{displaymath}
\mathbb P(\Omega(m)^c)
\leqslant\frac{\mathfrak c_{\Omega}}{(NT)^{6 +\kappa}}.
\end{displaymath}
On the other hand, as in the proof of Comte and Genon-Catalot \cite{CGC18}, Lemma 5, under Assumption \ref{assumption_basis_2} and on $\widehat\Lambda_{\kappa}(m)^c$,
\begin{displaymath}
L(m)\|\mathbf\Psi(m)^{-1}\|_{\textrm{op}}
\leqslant
\frac{\mathfrak c_{\kappa,T}}{2}\cdot\frac{NT}{\log(NT)}
\end{displaymath}
and
\begin{displaymath}
L(m)\|\widehat{\mathbf\Psi}(m)^{-1}\|_{\textrm{op}} >
\mathfrak c_{\kappa,T}\frac{NT}{\log(NT)},
\end{displaymath}
and then
\begin{eqnarray*}
\widehat\Lambda_{\kappa}(m)^c
 & \subset & \left\{
 L(m)\|\widehat{\mathbf\Psi}(m)^{-1} -\mathbf\Psi(m)^{-1}\|_{\textrm{op}}
 >\frac{\mathfrak c_{\kappa,T}}{2}\cdot
 \frac{NT}{\log(NT)}\right\}\\
 & \subset &
 \{\|\widehat{\mathbf\Psi}(m)^{-1} -\mathbf\Psi(m)^{-1}\|_{\textrm{op}} >
 \|\mathbf\Psi(m)^{-1}\|_{\textrm{op}}\}.
\end{eqnarray*}
Finally, by Comte and Genon-Catalot \cite{CGC18}, Proposition 4.(ii) which remains true for $H > 1/2$ without additional arguments,
\begin{displaymath}
\{\|\widehat{\mathbf\Psi}(m)^{-1} -\mathbf\Psi(m)^{-1}\|_{\textrm{op}} >
\|\mathbf\Psi(m)^{-1}\|_{\textrm{op}}\}
\subset\Omega(m)^c.
\end{displaymath}
Therefore,
\begin{displaymath}
\mathbb P(\widehat\Lambda_{\kappa}(m)^c)
\leqslant\mathbb P(\Omega(m)^c)
\leqslant\frac{\mathfrak c_{\Omega}}{(NT)^{6 +\kappa}}.
\end{displaymath}
%


%
\subsection{Proof of Corollary \ref{risk_bound}}
First of all,
\begin{equation}\label{risk_bound_1}
\mathbb E(\|\widetilde b_m - b_A\|_{f_T}^{2}) =
\mathfrak m_1 +\mathfrak m_2 +\mathfrak m_3
\end{equation}
where
\begin{eqnarray*}
 \mathfrak m_1 & := &
 \mathbb E(\|b_A\|_{f_T}^{2}\mathbf 1_{\widehat\Lambda_{\kappa}(m)^c}),\\
 \mathfrak m_2 & := &
 \mathbb E(\|\widehat b_m - b_A\|_{f_T}^{2}\mathbf 1_{\widehat\Lambda_{\kappa}(m)\cap\Omega(m)})
 \textrm{ and}\\
 \mathfrak m_3 & := &
 \mathbb E(\|\widehat b_m - b_A\|_{f_T}^{2}\mathbf 1_{\widehat\Lambda_{\kappa}(m)\cap\Omega(m)^c}).
\end{eqnarray*}
On the one hand, by Lemma \ref{bound_Lambda}, there exists a deterministic constant $\mathfrak c_1 > 0$, not depending on $m$, $N$ and $T$, such that
\begin{equation}\label{risk_bound_2}
\mathfrak m_1\leqslant\frac{\mathfrak c_1}{NT},
\end{equation}
and as in the proof of Comte and Genon-Catalot \cite{CGC19}, Proposition 1, Theorem \ref{empirical_risk_bound} allows to get
\begin{equation}\label{risk_bound_3}
\mathfrak m_2
\leqslant
\left(1 +\frac{4\mathfrak c_{\kappa,T}}{\log(NT)}\right)
\inf_{\tau\in\mathcal S_m}\|\tau - b_A\|_{f_T}^{2} +\frac{8}{NT}\textrm{trace}(\mathbf\Psi(m)^{-1}\mathbf\Psi(m,\sigma)).
\end{equation}
 On the other hand, by Lemma \ref{bound_e},
\begin{eqnarray*}
 \mathbb E(|\widehat{\mathbf x}(m)^*\widehat{\mathbf x}(m)|^2)
 & \leqslant &
 4m\sum_{j = 0}^{m - 1}\mathbb E(\langle\varphi_j,b\rangle_{N}^{4})+
 4\mathbb E(|\mathbf e(m)^*\mathbf e(m)|^2)\\
 & \leqslant &
 4mL(m)^2\int_{-\infty}^{\infty}b_A(x)^4f_T(x)dx\\
 & & +
 4\mathfrak c_{\mathbf e}\sigma^4
 \frac{\overline{\mathfrak m}_{4,H,M}(T)}{N^2T^{4 - 4H}}m(L(m)^2 + R(m)^2)\\
 & \leqslant &
 \mathfrak c_2
 \left(mL(m)^2 +
 \frac{\overline{\mathfrak m}_{4,H,M}(T)}{N^2T^{4 - 4H}}m(L(m)^2 + \rho^2L(m)^{2\kappa})\right)
\end{eqnarray*}
where
\begin{displaymath}
\mathfrak c_2 :=
\left(4\int_{-\infty}^{\infty}b_A(x)^4f_T(x)dx\right)
\vee
(4\mathfrak c_{\mathbf e}\sigma^4).
\end{displaymath}
As in the proof of Comte and Genon-Catalot \cite{CGC19}, Proposition 1, $\|\mathbf\Psi(m)\|_{\textrm{op}}\leqslant L(m)$, and by the definition of $\widehat\Lambda_{\kappa}(m)$,
\begin{eqnarray*}
 \mathbb E(\|\widehat b_m\|_{f_T}^{2}\mathbf 1_{\widehat\Lambda_{\kappa}(m)\cap\Omega(m)^c})
 & \leqslant &
 \|\mathbf\Psi(m)\|_{\textrm{op}}
 \mathbb E(
 \|\widehat{\mathbf\Psi}(m)^{-1}\|_{\textrm{op}}^{2}
 \widehat{\mathbf x}(m)^*\widehat{\mathbf x}(m)
 \mathbf 1_{\widehat\Lambda_{\kappa}(m)\cap\Omega(m)^c})\\
 & \leqslant &
 \frac{\mathfrak c_{\kappa,T}^{2}}{L(m)}
 \left(\frac{NT}{\log(NT)}\right)^2
 \mathbb P(\Omega(m)^c)^{1/2}
 \mathbb E(|\widehat{\mathbf x}(m)^*\widehat{\mathbf x}(m)|^2)^{1/2}\\
 & \leqslant &
 \frac{\mathfrak c_{2}^{1/2}\mathfrak c_{\kappa,T}^{2}}{L(m)}
 \left(\frac{NT}{\log(NT)}\right)^2
 \frac{\mathfrak c_{\Omega}^{1/2}}{(NT)^{3 +\kappa/2}}\\
 & &
 \times\left(m^{1/2}L(m) +
 \frac{\overline{\mathfrak m}_{2,H,M}(T)}{NT^{2 - 2H}}m^{1/2}(L(m) +\rho L(m)^{\kappa})\right).
\end{eqnarray*}
Then, there exists a deterministic constant $\mathfrak c_3 > 0$, not depending on $m$, $N$ and $T$, such that
\begin{equation}\label{risk_bound_4}
\mathfrak m_3
\leqslant 2\mathbb E(\|\widehat b_m\|_{f_T}^{2}\mathbf 1_{\widehat\Lambda_{\kappa}(m)\cap\Omega(m)^c}) +
2\|b_A\|_{f_T}^{2}\mathbb P(\Omega(m)^c)\\
\leqslant
\mathfrak c_3
\frac{\overline{\mathfrak m}_{2,H,M}(T)}{NT}.
\end{equation}
Inequalities (\ref{risk_bound_2}), (\ref{risk_bound_3}) and (\ref{risk_bound_4}) together with Inequality (\ref{risk_bound_1}) allow to conclude.
%


%
\subsection{Proof of Proposition \ref{bound_trace_term}}
First note that $\mathbf \Psi(m)$ and $\mathbf \Psi(m, \sigma)$ are symmetric and nonnegative, as for any vector $y$, $y^*\mathbf \Psi(m)y\geqslant 0$ and $y^*\mathbf \Psi(m, \sigma)y\geqslant 0$. 
Indeed for $y\in\mathbb R^m$, $y^*\mathbf\Psi_m y= \int \tau^2(x)f_T(x)dx\geqslant 0$ and 
\begin{eqnarray*}
 y^*\mathbf\Psi(m,\sigma) y 
 & = &
 \frac{\sigma^2}{T}
 \sum_{j,k = 0}^{m - 1}y^jy^k
 \mathbb E\left(\left(\int_{0}^{T}\varphi_j(X(s))\delta B(s)\right)
 \left(\int_{0}^{T}\varphi_k(X(s))\delta B(s)\right)\right)\\
 & = & \frac{\sigma^2}T {\mathbb E}\left[\left(\int_0^T\tau(X(s))\delta B(s)\right)^2\right]\geqslant 0, \quad
 \textrm{ where }
 \tau :=\sum_{j = 0}^{m - 1}y^j\varphi_j.
\end{eqnarray*}
Therefore, 
\begin{eqnarray*}
 \textrm{trace}(\mathbf\Psi(m)^{-1}\mathbf\Psi(m,\sigma))
 & \leqslant & \|\mathbf \Psi(m)^{-1}\|_{\textrm{op}}  \textrm{trace}(\mathbf\Psi(m,\sigma)) \\ 
 &=& \frac{\sigma^2}T  \|\mathbf \Psi(m)^{-1}\|_{\textrm{op}}  \sum_{j=0}^{m-1} {\mathbb E}\left[ \left(\int_0^T \varphi_j(X(s))\delta B(s)\right)^2\right]
\end{eqnarray*}
 By Theorem \ref{control_divergence_integral}, Equality (\ref{bound_e_1}) and Inequality (\ref{bound_e_2}),
\begin{eqnarray}
 \mathbb E\left[\left(\int_{0}^{T}\varphi_j(X(s))\delta B(s)\right)^2\right] 
 & \leqslant &
 \mathfrak c_{2,H,\sigma}
 \overline{\mathfrak m}_{2,H,M}(T)
 \left[\left(\int_{0}^{T}\mathbb E(|\varphi_j(X(s))|^{1/H})ds\right)^{2H}\right.
 \nonumber\\
 & &
 \left.+\left(\int_{0}^{T}\mathbb E(|\varphi_j'(X(s))|^2)^{1/(2H)}ds\right)^{2H}\right]
 \nonumber\\
 \label{bound_trace_term_2}
 &\leqslant &
 \mathfrak c_{2,H,\sigma}\overline{\mathfrak m}_{2,H,M}(T)T^{2H}
 \left(\int_{-\infty}^{\infty}\varphi_j(x)^2f_T(x)dx
 +\int_{-\infty}^{\infty}\varphi_j'(x)^2f_T(x)dx\right).
\end{eqnarray}
Finally we get
\begin{displaymath}
\textrm{trace}(\mathbf\Psi(m)^{-1}\mathbf\Psi(m,\sigma))\leqslant
\mathfrak c_{2,H,\sigma}\overline{\mathfrak m}_{2,H,M}(T)\frac{\sigma^2}{ T^{1-2H}}\|\mathbf\Psi(m)^{-1}\|_{\textrm{op}}(L(m)+ R(m)).
\end{displaymath}
\noindent For the second bound, let us denote by $\mathbf\Psi(m)^{-1/2}$ a symmetric square root of $\mathbf\Psi(m)^{-1}$, and write
\begin{eqnarray}
 \textrm{trace}(\mathbf\Psi(m)^{-1}\mathbf\Psi(m,\sigma))
 & = & \textrm{trace}(\mathbf\Psi(m)^{-1/2}\mathbf\Psi(m,\sigma)\mathbf\Psi(m)^{-1/2})
 \nonumber\\
 & \leqslant & m\|\mathbf\Psi(m)^{-1/2}\mathbf\Psi(m,\sigma)\mathbf\Psi(m)^{-1/2}\|_{{\rm op}}
 \nonumber\\
 \label{bound_trace_term_1}
 & = & m\sup_{\|x\|_{2,m} = 1}
 x^*\mathbf\Psi(m)^{-1/2}\mathbf\Psi(m,\sigma)\mathbf\Psi(m)^{-1/2}x
 = m\sup_{\|\mathbf\Psi(m)^{1/2}y\|_{2,m} = 1}
 y^*\mathbf\Psi(m,\sigma)y.
\end{eqnarray}
As already noticed, for any $y\in\mathbb R^m$ such that $\|\mathbf\Psi(m)^{1/2}y\|_{2,m} = 1$,
\begin{eqnarray*}
 y^*\mathbf\Psi(m,\sigma) y 
 & = & \frac{\sigma^2}T {\mathbb E}\left[\left(\int_0^T\tau(X(s))\delta B(s)\right)^2\right]
 \textrm{ with }
 \tau :=\sum_{j = 0}^{m - 1}y^j\varphi_j.
\end{eqnarray*}
By Theorem \ref{control_divergence_integral}, Equality (\ref{bound_e_1}) and Inequality (\ref{bound_e_2}) again,
\begin{equation}
 \mathbb E\left[\left(\int_{0}^{T}\tau(X(s))\delta B(s)\right)^2\right] 
 \label{bound_trace_term_2}
 \leqslant 
 \mathfrak c_{2,H,\sigma}\overline{\mathfrak m}_{2,H,M}(T)T^{2H}
 \left(\int_{-\infty}^{\infty}\tau(x)^2f_T(x)dx
 +\int_{-\infty}^{\infty}\tau'(x)^2f_T(x)dx\right).
\end{equation}
Since
\begin{displaymath}
\int_{-\infty}^{\infty}
\tau(x)^2f_T(x)dx =
y^*\mathbf\Psi(m)y = 1,
\end{displaymath}
by Inequalities (\ref{bound_trace_term_1}) and (\ref{bound_trace_term_2}),
\begin{displaymath}
\textrm{trace}(\mathbf\Psi(m)^{-1}\mathbf\Psi(m,\sigma))
\leqslant
\mathfrak c_{2,H,\sigma}\sigma^2
\frac{\overline{\mathfrak m}_{2,H,M}(T)}{T^{1 - 2H}}
m(1 + \|\mathbf\Psi(m)^{-1/2}
\Psi^{\ast}(m,\sigma)\mathbf\Psi(m)^{-1/2}\|_{{\rm op}})
\end{displaymath}
where
\begin{displaymath}
\Psi^{\ast}(m,\sigma) :=\left(\int_{-\infty}^{\infty}
\varphi_j'(x)\varphi_k'(x)f_T(x)dx\right)_{j,k = 0,\dots,m - 1}.
\end{displaymath}
Moreover, 
\begin{displaymath}
\|\mathbf\Psi(m)^{-1/2}\Psi^\ast(m,\sigma)\mathbf\Psi(m)^{-1/2}\|_{\textrm{op}}
\leqslant
\|\mathbf\Psi(m)^{-1/2}\|_{\textrm{op}}^{2}R(m).
\end{displaymath}
Therefore,
\begin{displaymath}
\textrm{trace}(\mathbf\Psi(m)^{-1}\mathbf\Psi(m,\sigma))
\leqslant
\mathfrak c_{2,H,\sigma}\sigma^2\frac{\overline{\mathfrak m}_{2,H,M}(T)}{T^{1 - 2H}}
m(1 +\|\mathbf \Psi(m)^{-1}\|_{\textrm{op}} R(m)).
\end{displaymath}
This concludes the proof.
%


%
\subsection{Proof of Proposition \ref{approximation_Skorokhod}}
Consider $x\in\mathbb R$ and $\varepsilon,t > 0$. For every $s\in [0,t]$,
\begin{displaymath}
\partial_x X_x(s) = 1 +
\int_{0}^{s}b'(X_x(r))\partial_x X_x(r)dr
\end{displaymath}
and, by Taylor's formula,
\begin{small}
\begin{displaymath}
 X_{x +\varepsilon}(s) - X_x(s) =\varepsilon
 +\int_{0}^{s}
 (X_{x +\varepsilon}(r) - X_x(r))
 \int_{0}^{1}b'(X_x(r) +\theta (X_{x +\varepsilon}(r) - X_x(r)))d\theta dr.
\end{displaymath}
\end{small}
\newline
So, for every $(u,v)\in [0,t]^2$ such that $v < u$,
\begin{displaymath}
\frac{\partial_x X_x(u)}{\partial_x X_x(v)} =
\exp\left(\int_{v}^{u}b'(X_x(r))dr\right)
\end{displaymath}
and
\begin{equation}\label{approximation_Skorokhod_1}
\frac{X_{x +\varepsilon}(u) - X_x(u)}{X_{x +\varepsilon}(v) - X_x(v)} =
\exp\left(\int_{v}^{u}\int_{0}^{1}b'(X_x(r) +\theta(X_{x +\varepsilon}(r) - X_x(r)))d\theta dr\right).
\end{equation}
For a given $\varphi\in\textrm{Lip}_{b}^{1}(\mathbb R)$, by Comte and Marie \cite{CM19}, Proposition 2.7,
\begin{displaymath}
 \Delta_{\varphi}^{S}(x,\varepsilon,t)\leqslant
 \alpha_H\sigma^2\int_{0}^{t}\int_{0}^{u}
 |\varphi'(X_x(u))|\Delta_{\varphi}(x,\varepsilon,u,v)(u - v)^{2H - 2}dvdu,
\end{displaymath}
where
\begin{displaymath}
\Delta_{\varphi}^{S}(x,\varepsilon,t) :=
\left|\int_{0}^{t}\varphi(X_x(u))\delta X_x(u) - S_{\varphi}(x,\varepsilon,t)\right|
\end{displaymath}
and, for every $(u,v)\in [0,t]^2$ such that $v < u$,
\begin{displaymath}
\Delta_{\varphi}(x,\varepsilon,u,v) :=
\left|
\frac{\partial_x X_x(u)}{\partial_x X_x(v)} -
\frac{X_{x +\varepsilon}(u) - X_x(u)}{X_{x +\varepsilon}(v) - X_x(v)}\right|.
\end{displaymath}
Since $b'(\mathbb R)\subset [m,M]$ and $b$ is two times continuously differentiable,
\begin{eqnarray*}
 \Delta_{\varphi}(x,\varepsilon,u,v) & = &
 \left|\exp\left(\int_{v}^{u}b'(X_x(r))dr\right)\right.\\
 & &
 -\left.
 \exp\left(\int_{v}^{u}\int_{0}^{1}b'(X_x(r) +\theta(X_{x +\varepsilon}(r) - X_x(r)))d\theta dr\right)\right|\\
 & \leqslant &
 \sup_{z\in [m(u - v),M(u - v)]}
 e^{z}\\
 & &
 \times
 \int_{v}^{u}\left|b'(X_x(r)) -
 \int_{0}^{1}b'(X_x(r) +\theta(X_{x +\varepsilon}(r) - X_x(r)))d\theta\right|dr\\
 & \leqslant &
 (1\vee e^{Mt})\int_{v}^{u}\int_{0}^{1}|b'(X_x(r)) -
 b'(X_x(r) +\theta(X_{x +\varepsilon}(r) - X_x(r)))|d\theta dr\\
 & \leqslant &
 \frac{\|b''\|_{\infty}}{2}
 (1\vee e^{Mt})\int_{v}^{u}|X_{x +\varepsilon}(r) - X_x(r)|dr.
\end{eqnarray*}
Consider $s\in\mathbb R_+$. By Equation (\ref{main_equation}),
\begin{eqnarray*}
 (X_{x +\varepsilon}(s) - X_x(s))^2 & = &
 \varepsilon^2 + 2\int_{0}^{s}
 (X_{x +\varepsilon}(r) - X_x(r))d(X_{x +\varepsilon} - X_x)(r)\\
 & = &
 \varepsilon^2 + 2\int_{0}^{s}
 (X_{x +\varepsilon}(r) - X_x(r))(b(X_{x +\varepsilon}(r)) - b(X_x(r)))dr.
\end{eqnarray*}
By the mean value theorem, there exists $x_s\in\mathbb R$ such that
\begin{eqnarray*}
 \partial_s(X_{x +\varepsilon}(s) - X_x(s))^2 & = &
 2(X_{x +\varepsilon}(s) - X_x(s))^2
 \frac{b(X_{x +\varepsilon}(s)) - b(X_x(s))}{X_{x +\varepsilon}(s) - X_x(s)}\\
 & = &
 2(X_{x +\varepsilon}(s) - X_x(s))^2
 b'(x_s)\leqslant
 2M(X_{x +\varepsilon}(s) - X_x(s))^2
\end{eqnarray*}
and then,
\begin{equation}\label{approximation_Skorokhod_2}
|X_{x +\varepsilon}(s) - X_x(s)|
\leqslant
\varepsilon e^{Ms}.
\end{equation}
 Therefore,
\begin{eqnarray*}
 \Delta_{\varphi}(x,\varepsilon,u,v)
 & \leqslant &
 \frac{\|b''\|_{\infty}}{2}\varepsilon
 (1\vee e^{Mt})
 \int_{v}^{u}e^{Mr}dr\\
 & \leqslant &
 \frac{\|b''\|_{\infty}}{2}\varepsilon(1\vee e^{Mt})\mathfrak m_M(u,v)
\end{eqnarray*}
where
\begin{displaymath}
\mathfrak m_M(u,v) :=
-\frac{1}{M}e^{Mv}\mathbf 1_{M < 0} +
(u - v)\mathbf 1_{M = 0} +
\frac{1}{M}e^{Mu}\mathbf 1_{M > 0}.
\end{displaymath}
Finally, using the above bounds, and in a second stage, the integration by parts formula, we get
\begin{eqnarray*}
 \Delta_{\varphi}^{S}(x,\varepsilon,t) & \leqslant &
 \alpha_H\sigma^2\frac{\|b''\|_{\infty}}{2}\varepsilon(1\vee e^{Mt})
 \int_{0}^{t}\int_{0}^{u}|\varphi'(X_x(u))|\mathfrak m_M(u,v)(u - v)^{2H - 2}dvdu\\
 & \leqslant &
 \alpha_H\sigma^2\frac{\|b''\|_{\infty}\|\varphi'\|_{\infty}}{2}\varepsilon(1\vee e^{Mt})
 \left(-\frac{1}{M}\mathbf 1_{M < 0}
 \int_{0}^{t}e^{Mv}\int_{v}^{t}(u - v)^{2H - 2}dudv\right.\\
 & &
 \left. +
 \mathbf 1_{M = 0}
 \int_{0}^{t}\int_{0}^{u}(u - v)^{2H - 2}dvdu +
 \frac{1}{M}\mathbf 1_{M > 0}
 \int_{0}^{t}e^{Mu}\int_{0}^{u}(u - v)^{2H - 2}dvdu
 \right)\\
 & \leqslant &
 \alpha_H\sigma^2\frac{\|b''\|_{\infty}\|\varphi'\|_{\infty}}{2}\varepsilon(1\vee e^{Mt})
 \left(\frac{1}{M^2}\mathbf 1_{M < 0}
 \left(\frac{t^{2H - 1}}{2H - 1} -\int_{0}^{t}e^{Mv}(t - v)^{2H - 2}dv\right)\right.\\
 & &
 \left. +
 \frac{1}{2H}\mathbf 1_{M = 0}
 \int_{0}^{t}u^{2H}du +
 \frac{1}{M^2}\mathbf 1_{M > 0}
 \left(\frac{e^{Mt}t^{2H - 1}}{2H - 1} -\int_{0}^{t}e^{Mu}u^{2H - 2}du\right)
 \right)\\
 & \leqslant &
 \alpha_H\sigma^2\frac{\|b''\|_{\infty}\|\varphi'\|_{\infty}}{2}\varepsilon t^{2H - 1}\mathfrak m_{H,M}(t)
\end{eqnarray*}
where
\begin{displaymath}
\mathfrak m_{H,M}(t) =
\frac{1}{M^2(2H - 1)}\mathbf 1_{M < 0} +
\frac{t^2}{2H(2H + 1)}\mathbf 1_{M = 0} +
\frac{e^{2Mt}}{M^2(2H - 1)}\mathbf 1_{M > 0}.
\end{displaymath}
%


%
\subsection{Proof of Corollary \ref{risk_bound_computable_estimator}}
First of all,
\begin{displaymath}
\mathbb E(\|\widetilde b_{m,\varepsilon} - b_A\|_{f_T}^{2})
\leqslant
2\mathbb E(\|\widehat b_{m,\varepsilon} -\widehat b_m\|_{f_T}^{2}\mathbf 1_{\widehat\Lambda_{\kappa}(m)}) +
2\mathbb E(\|\widetilde b_m - b_A\|_{f_T}^{2})
\end{displaymath}
and
\begin{eqnarray*}
 \mathbb E(\|\widetilde b_m - b_A\|_{f_T}^{2})
 & \leqslant &
 \left(1 + 4\frac{3\log(3/2)-1}{(7 +\kappa)\log(NT)}\right)
 \inf_{\tau\in\mathcal S_m}
 \|\tau - b_A\|_{f_T}^{2}\\
 & &
 +\frac{8}{NT}\normalfont{\textrm{trace}}(\mathbf\Psi(m)^{-1}\mathbf\Psi(m,\sigma)) +
 \overline{\mathfrak c}_{\rho,\kappa,\sigma}(1 +\mathfrak b_T)\frac{\overline{\mathfrak m}_{2,H,M}(T)}{NT}
\end{eqnarray*}
by Corollary \ref{risk_bound}. So, let us find a suitable control for $\mathbb E(\|\widehat b_{m,\varepsilon} -\widehat b_m\|_{f_T}^{2}\mathbf 1_{\widehat\Lambda_{\kappa}(m)})$. On the one hand,
\begin{eqnarray*}
 \|\widehat b_{m,\varepsilon} -\widehat b_m\|_{f_T}^{2}
 & = &
 \|\langle\widehat\theta(m,\varepsilon) -\widehat\theta(m),(\varphi_0,\dots,\varphi_{m - 1})(.)\rangle_{2,m}\|_{f_T}^{2}\\
 & \leqslant &
 \|\widehat\theta(m,\varepsilon) -\widehat\theta(m)\|_{2,m}^{2}
 \int_A\|(\varphi_0,\dots,\varphi_{m - 1})(x)\|_{2,m}^{2}f_T(x)dx\\
 & \leqslant &
 \|\widehat\theta(m,\varepsilon) -\widehat\theta(m)\|_{2,m}^{2}L(m).
\end{eqnarray*}
 On the other hand, by Proposition \ref{approximation_Skorokhod}, for every $i\in\{1,\dots,N\}$ and $j\in\{0,\dots,m - 1\}$,
\begin{displaymath}
\left|S_{\varphi_j}^{i}(x_0,\varepsilon,T) -\int_{0}^{T}\varphi_j(X_{x_0}^{i}(s))\delta X_{x_0}^{i}(s)\right|
\leqslant
\alpha_H\sigma^2\frac{\|b''\|_{\infty}\|\varphi'\|_{\infty}}{2}\varepsilon T^{2H - 1}\mathfrak m_{H,M}(T).
\end{displaymath}
Then,
\begin{eqnarray*}
 \|\widehat\theta(m,\varepsilon) -\widehat\theta(m)\|_{2,m}^{2}
 & \leqslant &
 \|\widehat{\mathbf\Psi}(m)^{-1}\|_{\textrm{op}}^{2}
 \|\widehat{\mathbf x}(m,\varepsilon) -\widehat{\mathbf x}(m)\|_{2,m}^{2}\\
 & = &
 \|\widehat{\mathbf\Psi}(m)^{-1}\|_{\textrm{op}}^{2}
 \sum_{j = 0}^{m - 1}
 \left|\frac{1}{NT}\sum_{i = 1}^{N}\left(S_{\varphi_j}^{i}(x_0,\varepsilon,T) -\int_{0}^{T}\varphi_j(X_{x_0}^{i}(s))\delta X_{x_0}^{i}(s)\right)\right|^2\\
 & \leqslant &
 \mathfrak c_1\|\widehat{\mathbf\Psi}(m)^{-1}\|_{\textrm{op}}^{2}\varepsilon^2T^{4H - 4}\mathfrak m_{H,M}(T)^2R(m)
\end{eqnarray*}
with
\begin{displaymath}
\mathfrak c_1 :=
\alpha_{H}^{2}\sigma^4\frac{\|b''\|_{\infty}^{2}}{4}.
\end{displaymath}
Therefore, by the definition of $\widehat\Lambda_{\kappa}(m)$, there exists a deterministic constant $\mathfrak c_{\rho,\kappa,\sigma,H,\mathfrak b''} > 0$, depending only on $\rho$, $\kappa$, $\sigma$, $H$ and $\mathfrak b''$, such that
\begin{eqnarray*}
 2\mathbb E(\|\widehat b_{m,\varepsilon} -\widehat b_m\|_{f_T}^{2}\mathbf 1_{\widehat\Lambda_{\kappa}(m)})
 & \leqslant &
 2\mathfrak c_1\mathfrak c_{\kappa,T}^{2}\left(\frac{NT}{\log(NT)}\right)^2L(m)^{-1}\varepsilon^2T^{4H - 4}\mathfrak m_{H,M}(T)^2R(m)\\
 & \leqslant &
 \mathfrak c_{\rho,\kappa,\sigma,H,\mathfrak b''}
 L(m)^{\kappa -1}
 \left(\frac{NT}{\log(NT)}\right)^2\frac{\mathfrak m_{H,M}(T)^2}{T^{4 - 4H}}\varepsilon^2.
\end{eqnarray*}
%


%
\subsection{Proof of Proposition \ref{empirical_risk_bound_estimator_dag}}
First, let us prove (\ref{bounddag1}). We have
\begin{displaymath}
\|\widetilde b_{m}^{\dag,\prime} - b_A'\|_{N}^{2} =
\|\widehat b_{m}^{\dag,\prime} - b_A'\|_{N}^{2}\mathbf 1_{\widehat\Lambda_0(m)} +
\|b_A'\|_{N}^{2}\mathbf 1_{\widehat\Lambda_0(m)^c}.
\end{displaymath}
On the one hand, by Lemma \ref{bound_Lambda},
\begin{displaymath}
\mathbb E(\|b_A'\|_{N}^{2}\mathbf 1_{\widehat\Lambda_0(m)^c})
\leqslant
\mathbb E(\|b_A'\|_{N}^{4})^{1/2}
\mathbb P(\widehat\Lambda_0(m)^c)^{1/2}
\leqslant
\frac{\mathfrak c_1}{NT}
\end{displaymath}
where
\begin{displaymath}
\mathfrak c_1 :=
\mathfrak c_{\Omega}^{1/2}
\left(\int_{-\infty}^{\infty}b_A'(x)^4f_T(x)dx\right)^{1/2}
\leqslant
\mathfrak c_{\Omega}^{1/2}\|b'\|_{\infty}^{2} =
\mathfrak c_{\ref{empirical_risk_bound_estimator_dag}}.
\end{displaymath}
On the other hand, since Model (\ref{alternative_model}) is noiseless, $\widehat b_{m}^{\dag,\prime} = p_{N,\mathcal S_m}^{\perp}(b_A')$, and then
\begin{displaymath}
\|\widehat b_{m}^{\dag,\prime} - b_A'\|_{N}^{2} =
\inf_{\tau\in\mathcal S_m}\|\tau - b_A'\|_{N}^{2}.
\end{displaymath}
Therefore,
\begin{displaymath}
\mathbb E(\|\widetilde b_{m}^{\dag,\prime} - b_A'\|_{N}^{2})\leqslant
\inf_{\tau\in\mathcal S_m}\|\tau - b_A'\|_{f_T}^{2} +
\frac{\mathfrak c_{\ref{empirical_risk_bound_estimator_dag}}}{NT}.
\end{displaymath}
Now let us prove (\ref{bounddag2}). We start by writing that 
\begin{equation}\label{empirical_risk_bound_estimator_dag_1}
\mathbb E(\|\widetilde b_{m}^{\dag,\prime} - b_A'\|_{f_T}^{2}) =
\mathfrak m_1 +\mathfrak m_2 +\mathfrak m_3
\end{equation}
where
\begin{eqnarray*}
 \mathfrak m_1 & := &
 \mathbb E(\|b_A'\|_{f_T}^{2}\mathbf 1_{\widehat\Lambda_0(m)^c}),\\
 \mathfrak m_2 & := &
 \mathbb E(\|\widehat b_{m}^{\dag,\prime} - b_A'\|_{f_T}^{2}\mathbf 1_{\widehat\Lambda_0(m)\cap\Omega(m)})
 \textrm{ and}\\
 \mathfrak m_3 & := &
 \mathbb E(\|\widehat b_{m}^{\dag,\prime} - b_A'\|_{f_T}^{2}\mathbf 1_{\widehat\Lambda_0(m)\cap\Omega(m)^c}).
\end{eqnarray*}
Let us find suitable bounds on $\mathfrak m_1$, $\mathfrak m_2$ and $\mathfrak m_3$.
\begin{itemize}
 \item\textbf{Bound on $\mathfrak m_1$.} By Lemma \ref{bound_Lambda}, there exists a deterministic constant $\mathfrak c_1 > 0$, not depending on $m$, $N$ and $T$, such that
 \begin{equation}\label{empirical_risk_bound_estimator_dag_2}
 \mathfrak m_1\leqslant\frac{\mathfrak c_1}{NT}.
 \end{equation}
 \item\textbf{Bound on $\mathfrak m_2$.} By denoting $p_{f_T,\mathcal S_m}^{\perp}$ the orthogonal projection from $\mathbb L^2(A,f_T(x)dx)$ onto $\mathcal S_m$ with respect to the $f_T$-weighted scalar product $\langle .,.\rangle_{f_T}$,
 \begin{eqnarray*}
  \|\widehat b_{m}^{\dag,\prime} - b_A'\|_{f_T}^{2} & = &
  \|b_A' - p_{N,\mathcal S_m}^{\perp}(b_A')\|_{f_T}^{2} =
  \|b_A' - p_{f_T,\mathcal S_m}^{\perp}(b_A') -
  p_{N,\mathcal S_m}^{\perp}(b_A' - p_{f_T,\mathcal S_m}^{\perp}(b_A'))\|_{f_T}^{2}\\
  & = &
  \inf_{\tau\in\mathcal S_m}\|\tau - b_A'\|_{f_T}^{2} +
  \|p_{N,\mathcal S_m}^{\perp}(g)\|_{f_T}^{2}
 \end{eqnarray*}
 with $g := b_A' - p_{f_T,\mathcal S_m}^{\perp}(b_A')$. By Comte and Genon-Catalot \cite{CGC19}, Lemma 5 applied to this function $g$,
 \begin{eqnarray*}
  \mathbb E(\|p_{N,\mathcal S_m}^{\perp}(g)\|_{f_T}^{2}
  \mathbf 1_{\widehat\Lambda_0(m)\cap\Omega(m)})
  & \leqslant &
  \frac{2\mathfrak c_{0,T}}{\log(NT)}\|g\|_{f_T}^{2} =
  \frac{2\mathfrak c_{0,T}}{\log(NT)}\inf_{\tau\in\mathcal S_m}\|\tau - b_A'\|_{f_T}^{2}.
 \end{eqnarray*}
 \item\textbf{Bound on $\mathfrak m_3$.} First of all, by Lemma \ref{bound_Lambda},
 \begin{displaymath}
 \mathfrak m_3
 \leqslant
 \frac{\mathfrak c_{\Omega}^{1/2}}{(NT)^3}
 \mathbb E(\|\widehat b_{m}^{\dag,\prime}\|_{f_T}^{4}\mathbf 1_{\widehat\Lambda_0(m)})^{1/2} +
 \|b_A'\|_{f_T}^{2}\frac{\mathfrak c_{\Omega}}{NT}.
 \end{displaymath}
 So, it remains to find a bound on $\mathbb E(\|\widehat b_{m}^{\dag,\prime}\|_{f_T}^{4}\mathbf 1_{\widehat\Lambda_0(m)})$. On the one hand,
 \begin{eqnarray*}
  \|\widehat b_{m}^{\dag,\prime}\|_{f_T}^{2}
  & = &
  \int_{-\infty}^{\infty}\left[\sum_{j = 0}^{m - 1}\widehat\theta_j'\varphi_j(x)\right]^2f_T(x)dx\\
  & = &
  \widehat\theta'(m)^{*}\mathbf\Psi(m)\widehat\theta'(m) =
  \widehat{\mathbf y}(m)^*\widehat{\mathbf\Psi}(m)^{-1}
  \mathbf\Psi(m)\widehat{\mathbf\Psi}(m)^{-1}\widehat{\mathbf y}(m) =
  \|\mathbf\Psi(m)^{1/2}\widehat{\mathbf\Psi}(m)^{-1}\widehat{\mathbf y}(m)\|_{2,m}^{2}\\
  & \leqslant &
  \|\mathbf\Psi(m)\|_{{\rm op}}
  \|\widehat{\mathbf\Psi}(m)^{-1}\|_{{\rm op}}^{2}
  \sum_{j = 0}^{m - 1}\langle\varphi_j,b'\rangle_{N}^{2}.
 \end{eqnarray*}
 On the other hand, $\|\mathbf\Psi(m)\|_{{\rm op}}\leqslant L(m)$,
 \begin{displaymath}
 \|\widehat{\mathbf\Psi}(m)^{-1}\|_{{\rm op}}^{2}
 \leqslant
 \mathfrak c_{0,T}^{2}L(m)^{-2}\left(\frac{NT}{\log(NT)}\right)^2
 \end{displaymath}
 on $\widehat{\Lambda}_0(m)$, and
 \begin{eqnarray*}
 & &\mathbb E\left[\left|\sum_{j = 0}^{m - 1}\langle\varphi_j,b'\rangle_{N}^{2}\right|^2\right]^{1/2}\\
 & &
 \hspace{1.5cm}\leqslant
 \sqrt m\mathbb E\left[
 \sum_{j = 0}^{m - 1}\left|\frac{1}{NT}\sum_{i = 1}^{N}\int_{0}^{T}\varphi_j(X_{x_0}^{i}(s))b'(X_{x_0}^{i}(s))ds\right|^4\right]^{1/2}\\
 & &
 \hspace{1.5cm}\leqslant
 \sqrt m\left[
 \frac{1}{T}\int_{0}^{T}\mathbb E\left[
 \sum_{j = 0}^{m - 1}
 \varphi_j(X_{x_0}^{1}(s))^4b'(X_{x_0}^{1}(s))^4\right]ds\right]^{1/2}\\
 & &
 \hspace{1.5cm}\leqslant
 \sqrt mL(m)\left(\frac{1}{T}\int_{0}^{T}\int_{-\infty}^{\infty}b'(x)^4p_s(x_0,x)dxds\right)^{1/2} =
 \sqrt mL(m)\left(\int_{-\infty}^{\infty}b'(x)^4f_T(x)dx\right)^{1/2}.
 \end{eqnarray*}
 Then,
 \begin{displaymath}
 \mathbb E(\|\widehat b_{m}^{\dag,\prime}\|_{f_T}^{4}\mathbf 1_{\widehat\Lambda_0(m)})^{1/2}
 \leqslant
 \mathfrak c_{0,T}\|b'\|_{\infty}^{2}\sqrt m\frac{NT}{\log(NT)}.
 \end{displaymath}
\end{itemize}
%


%
\subsection{Proof of Corollary \ref{risk_bound_computable_estimator_dag}}
First of all,
\begin{displaymath}
\mathbb E(\|\widetilde b_{m,\varepsilon}^{\dag,\prime} - b_A'\|_{f_T}^{2})
\leqslant
2\mathbb E(\|\widehat b_{m,\varepsilon}^{\dag,\prime} -\widehat b_{m}^{\dag,\prime}\|_{f_T}^{2}\mathbf 1_{\widehat\Lambda_0(m)}) +
2\mathbb E(\|\widetilde b_{m}^{\dag,\prime} - b_A'\|_{f_T}^{2})
\end{displaymath}
and
\begin{displaymath}
\mathbb E(\|\widetilde b_{m}^{\dag,\prime} - b_A'\|_{f_T}^{2})
\leqslant
\left(1 +  \frac{2T\mathfrak c_{0,T}}{\log(NT)}\right)
\inf_{\tau\in\mathcal S_m}
\|\tau - b_A'\|_{f_T}^{2}
+\overline{\mathfrak c}_{\mathfrak b'}\frac{1 +\mathfrak c_{0,T}}{NT}
\end{displaymath}
by (\ref{bounddag2}) of Proposition \ref{empirical_risk_bound_estimator_dag}. So, let us find a suitable control for $\mathbb E(\|\widehat b_{m,\varepsilon}^{\dag,\prime} -\widehat b_{m}^{\dag,\prime}\|_{f_T}^{2}\mathbf 1_{\widehat\Lambda_0(m)})$. On the one hand,
\begin{eqnarray*}
 \|\widehat b_{m,\varepsilon}^{\dag,\prime} -\widehat b_{m}^{\dag,\prime}\|_{f_T}^{2}
 & = &
 \|\langle\widehat\theta'(m,\varepsilon) -\widehat\theta'(m),(\varphi_0,\dots,\varphi_{m - 1})(.)\rangle_{2,m}\|_{f_T}^{2}\\
 & \leqslant &
 \|\widehat\theta'(m,\varepsilon) -\widehat\theta'(m)\|_{2,m}^{2}
 \int_A\|(\varphi_0,\dots,\varphi_{m - 1})(x)\|_{2,m}^{2}f_T(x)dx\\
 & \leqslant &
 \|\widehat\theta'(m,\varepsilon) -\widehat\theta'(m)\|_{2,m}^{2}L(m).
\end{eqnarray*}
On the other hand, for every $s\in [0,T]$, by (\ref{approximation_Skorokhod_1}) and (\ref{approximation_Skorokhod_2}),
\begin{eqnarray*}
 |\dot Y_{x_0}(s) -\dot Y_{x_0 +\varepsilon}(s)|
 & = &
 \left|b'(X_{x_0}(s)) -
 \int_{0}^{1}b'(X_{x_0}(s) +\theta(X_{x_0 +\varepsilon}(s) - X_{x_0}(s)))d\theta\right|\\
 & \leqslant &
 \int_{0}^{1}|b'(X_{x_0}(s)) -
 b'(X_{x_0}(s) +\theta(X_{x_0 +\varepsilon}(s) - X_{x_0}(s)))|d\theta\\
 & \leqslant &
 \|b''\|_{\infty}|X_{x_0 +\varepsilon}(s) - X_{x_0}(s)|\int_{0}^{1}\theta d\theta
 \leqslant
 \frac{\|b''\|_{\infty}}{2}e^{Ms}\varepsilon.
\end{eqnarray*}
Therefore, thanks to the definition of $\widehat\Lambda_0(m)$,
\begin{eqnarray*}
 & &
 \|\widehat\theta'(m,\varepsilon) -\widehat\theta'(m)\|_{2,m}^{2}\mathbf 1_{\widehat\Lambda_0(m)}\\
 & &
 \hspace{1.5cm}
 \leqslant
 \|\widehat{\mathbf\Psi}(m)^{-1}\|_{\textrm{op}}^{2}\mathbf 1_{\widehat\Lambda_0(m)}
 \sum_{j = 0}^{m - 1}\left|
 \frac{1}{NT}
 \sum_{i = 1}^{N}
 \int_{0}^{T}\varphi_j(X_{x_0}^{i}(s))d(Y_{x_0,\varepsilon}^{i} - Y_{x_0}^{i})(s)\right|^2\\
 & &
 \hspace{1.5cm}
 \leqslant
 \frac{\|b''\|_{\infty}^{2}e^{2MT}}{4}\varepsilon^2
 \|\widehat{\mathbf\Psi}(m)^{-1}\|_{\textrm{op}}^{2}\mathbf 1_{\widehat\Lambda_0(m)}
 \sum_{j = 0}^{m - 1}\left|
 \frac{1}{NT}
 \sum_{i = 1}^{N}
 \int_{0}^{T}\varphi_j(X_{x_0}^{i}(s))ds\right|^2\\
 & &
 \hspace{1.5cm}
 \leqslant
 \frac{\|b''\|_{\infty}^{2}e^{2MT}}{4}\varepsilon^2
 L(m)\|\widehat{\mathbf\Psi}(m)^{-1}\|_{\textrm{op}}^{2}\mathbf 1_{\widehat\Lambda_0(m)}
 \leqslant
 \frac{\|b''\|_{\infty}^{2}}{4}\mathfrak c_{0,T}^{2}e^{2MT}\frac{1}{L(m)}\left(\frac{NT}{\log(NT)}\right)^2\varepsilon^2.
\end{eqnarray*}
\textbf{Conflict of interest statement.} The authors certify that they have no affiliations with or involvement in any organization or entity with any financial interest, or non-financial interest in the subject matter or materials discussed in this manuscript.
\\
\\
\textbf{Acknowledgement.} We thank an anonymous referee for interesting remarks and helpful suggestion.
 

%

%
\end{document}